\documentclass[microtype]{gtpart}
\usepackage{pinlabel}
\usepackage{graphicx}
\usepackage{color}

\usepackage{amsfonts}
\usepackage{amssymb}
\usepackage{amsxtra}

\newtheorem{thm}{Theorem}[section]
\newtheorem{prop}[thm]{Proposition}
\newtheorem{lem}[thm]{Lemma}
\newtheorem{cor}[thm]{Corollary}
\newtheorem{conj}[thm]{Conjecture}

\theoremstyle{definition}
\newtheorem{defn}[thm]{Definition}
\newtheorem{rem}[thm]{Remark}

\newtheorem{exmp}[thm]{Example}
\newtheorem{ques}[thm]{Question}
\newtheorem{prob}[thm]{Problem}

\numberwithin{equation}{section}

\newcommand{\abs}[1]{\lvert{#1}\rvert}

\renewcommand{\bar}[1]{\overline{#1}}
\newcommand{\boundary}{\partial}
\newcommand{\gen}[1]{\langle{#1}\rangle}

\newcommand{\presentation}[2]{\langle\, {#1} \mid {#2} \,\rangle}

\newcommand{\set}[2]{\{\,{#1} \mid {#2} \,\}}
\newcommand{\bigset}[2]{ \bigl\{ \, {#1} \bigm| {#2} \, \bigr\} }
\renewcommand{\emptyset}{\varnothing}

\renewcommand{\setminus}{-}

\newcommand{\field}[1]{\mathbb{#1}}
\newcommand{\integers}{\Z}
\newcommand{\N}{\field{N}}

\DeclareMathOperator{\CAT}{CAT}

\newcommand{\ball}[2]{B ( {#1}, {#2} )}

\newcommand{\nbd}[2]{\mathcal{N}_{#2}({#1})}  
\newcommand{\bignbd}[2]{\mathcal{N}_{#2} \bigl( {#1} \bigr)}
\newcommand{\Set}[1]{\mathcal{#1}}

\DeclareMathOperator{\height}{height}
\DeclareMathOperator{\width}{width}

\DeclareMathOperator{\Cayley}{Cayley}
\DeclareMathOperator{\diam}{diam}
\DeclareMathOperator{\Sat}{Sat} 

\def\RomanianComma#1{\setbox0=\hbox{#1}{\ooalign{\hidewidth
    \lower1.2ex\hbox{$\mspace{1mu}^{\text{,}}$}\hidewidth\crcr\unhbox0}}}
\newcommand{\Drutu}{Dru{\RomanianComma{t}u}}

\usepackage{ifthen}

\newcommand{\showcomments}{yes}

\newsavebox{\commentbox}
{\ifthenelse{\equal{\showcomments}{yes}}%
{\footnotemark
    \begin{lrbox}{\commentbox}
    \begin{minipage}[t]{1.25in}\raggedright\sffamily\upshape\tiny
    \footnotemark[\arabic{footnote}]}
{\begin{lrbox}{\commentbox}}}%
{\ifthenelse{\equal{\showcomments}{yes}}%
{\end{minipage}\end{lrbox}\marginpar{\usebox{\commentbox}}}
{\end{lrbox}}}

\title{Packing subgroups in relatively hyperbolic groups}

\author[G\,C Hruska]{G Christopher Hruska}
\givenname{G Christopher}
\surname{Hruska}
\address{Department of Mathematical Sciences\\
University of Wisconsin--Milwaukee\\\newline
P.O.~Box 413\\
Milwaukee, WI 53201\\
USA}
\email{chruska@uwm.edu}
\urladdr{http://www.uwm.edu/~chruska}

\author[D\,T Wise]{Daniel T Wise}
\givenname{Daniel T}
\surname{Wise}
\address{Department of Mathematics and Statistics.\\
McGill University\\\newline
Montreal, Quebec, Canada H3A 2K6}
\email{wise@math.mcgill.ca}
\urladdr{http://www.math.mcgill.ca/~wise}

\keyword{relative hyperbolicity}
\keyword{quasiconvex subgroup}
\keyword{width}
\keyword{cube complex}
\subject{primary}{msc2000}{20F65}
\subject{secondary}{msc2000}{20F67}

\arxivreference{math.GR/0609369}
\arxivpassword{94fjd}

\begin{document}

\begin{abstract}
We introduce the bounded packing property for a subgroup of a countable
discrete group $G$.
This property gives a finite upper bound on
the number of left cosets of the subgroup that are pairwise close in $G$.
We establish basic properties of bounded packing,
and give many examples; for instance, every subgroup of a countable,
virtually nilpotent group has bounded packing.
We explain several natural connections between bounded
packing and group actions on $\CAT(0)$ cube complexes.

Our main result
establishes the bounded packing of relatively quasiconvex subgroups of a
relatively hyperbolic group, under mild hypotheses.
As an application, we prove that relatively quasiconvex subgroups
have finite height and width, properties that strongly restrict
the way families of distinct conjugates of the subgroup
can intersect.
We prove that an infinite, nonparabolic relatively quasiconvex
subgroup of a relatively hyperbolic group has finite index in its
commensurator.
We also prove a virtual malnormality theorem for separable, relatively quasiconvex subgroups, which is new even in the word hyperbolic case.
\end{abstract}

\begin{asciiabstract}
We introduce the bounded packing property for a subgroup of a countable
discrete group G.
This property gives a finite upper bound on
the number of left cosets of the subgroup that are pairwise close in G.
We establish basic properties of bounded packing,
and give many examples; for instance, every subgroup of a countable,
virtually nilpotent group has bounded packing.
We explain several natural connections between bounded
packing and group actions on CAT(0) cube complexes.

Our main result
establishes the bounded packing of relatively quasiconvex subgroups of a
relatively hyperbolic group, under mild hypotheses.
As an application, we prove that relatively quasiconvex subgroups
have finite height and width, properties that strongly restrict
the way families of distinct conjugates of the subgroup
can intersect.
We prove that an infinite, nonparabolic relatively quasiconvex
subgroup of a relatively hyperbolic group has finite index in its
commensurator.
We also prove a virtual malnormality theorem for separable, relatively quasiconvex subgroups, which is new even in the word hyperbolic case.
\end{asciiabstract}

\maketitle


\section{Introduction}
\label{sec:Introduction}

In treating groups as geometric objects,
it is natural to emphasize certain subgroups
for special treatment, and these ``quasiconvex'' subgroups
play a role akin to convex subspaces of a geodesic metric space.
While for arbitrary groups the quasiconvexity of a subgroup
is not even well-defined, for word-hyperbolic groups
quasiconvexity is independent of the choice of generating system,
and there is a useful theory of such subgroups which
endows them with the properties of the ambient group,
and then examines the favorable way in which they are embedded.

There has been substantial progress in the past few years
generalizing properties of word-hyperbolic groups
to properties of relatively hyperbolic groups
with hypotheses on the peripheral subgroups.
This paper is part of this trend
where we focus on generalizations of ideas related to
a certain important property of a quasiconvex subgroup.
 Let us first discuss the background of this property before
 indicating our generalizations.

An infinite quasiconvex subgroup $H$ of a word-hyperbolic group $G$
cannot be normal unless it is of finite index (see Short \cite{ABC91}).
More generally, there are substantial limits on the way it
intersects its conjugates. The definitive result in this direction,
obtained by Gitik--Mitra--Rips--Sageev,
states that a quasiconvex subgroup has ``finite width'' \cite{GMRS98}.
This roughly means that there is an upper bound on the number of conjugates
whose pairwise intersection is infinite.
More precisely the \emph{width} of $H$ in $G$ equals $(n-1)$
if $n$ is the smallest number with the property
that for any $n$ distinct right cosets $g_1H, \dots , g_nH$
the intersection $g_i H g_i^{-1} \cap g_j H g_j^{-1}$ is finite
for some $i,j$.
It is unknown whether having finite width characterizes
quasiconvexity, but it seems unlikely that one could jump from an
algebraic hypothesis to a geometric conclusion in this context.

\subsection{Bounded Packing}

The geometric mechanism lying behind the finite width conclusion
is what we call ``bounded packing.''
In a metric space $(X,d)$ the distance between two subsets $Y$ and $Z$
is the infimum of distances between points $y\in Y$ and $z\in Z$.
If $H$ is a subgroup of $G$, then $H$ has \emph{bounded packing in $G$} if for 
each
$D$, there is a bound $M=M(D)$ on the number of distinct cosets $g_iH$ such
that $d(g_iH, g_jH)<D$ for each $i,j$.
(A more precise statement can be found in Section~\ref{sec:BoundedPacking}
below.)
The connection between
width and the bounded packing of a $\kappa$-quasiconvex subgroup
$H$ of the $\delta$-hyperbolic group $G$, is that if
$g_iHg_i^{-1}\cap g_jHg_j^{-1}$ is infinite, then $d(g_iH,g_jH)<K=K(\kappa,\delta)$.

The goal of this paper is to prove a bounded packing statement
in the relatively hyperbolic context,
and to deduce from this an appropriate finite width consequence.

Bounded packing is the fundamental notion operating behind proofs about
the widths of subgroups and appears implicitly in
Gitik--Mitra--Rips--Sageev \cite{GMRS98}.
It is a natural algebraic generalization of the
the Finite Plane Intersection Property (FPIP) for
universal covers of surfaces in $3$--manifolds
(see, for instance, Rubinstein--Sageev \cite{RubinsteinSageev99}).
It is also the key point used by Sageev \cite{Sageev97}
to prove the finite
dimensionality of the cube complex arising from Sageev's construction
applied to a codimension--$1$ quasiconvex subgroup of a word-hyperbolic 
group, as discussed in Section~\ref{sec:Cubulating} below.
In \cite{HruskaWiseAxioms}, the authors
apply the bounded packing property
to generalize Sageev's finite dimensional cubulation
result to codimension--$1$ relatively quasiconvex subgroups of relatively hyperbolic groups.
Our main theorem about bounded packing is the following
result about relatively hyperbolic groups.

\begin{thm}\label{thm:RelHypBoundedPacking}
Let $H$ be a relatively quasiconvex subgroup of a relatively
hyperbolic group $G$.
Suppose $H\cap gPg^{-1}$ has bounded packing in $gPg^{-1}$ for each
conjugate of each peripheral subgroup $P$.
Then $H$ has bounded packing in $G$.
\end{thm}

This theorem generalizes a result of Gitik--Mitra--Rips--Sageev \cite{GMRS98}
on quasiconvex subgroups of a word hyperbolic group
and also generalizes a result of Rubinstein--Sageev \cite{RubinsteinSageev99}
on geometrically finite surfaces in $3$--manifolds whose toral decomposition
has only hyperbolic components.

We note that the condition about intersections with peripheral
subgroups is necessary.
Indeed, every group $G$ is hyperbolic relative to itself,
and in that case every subgroup of $G$ is relatively quasiconvex.

We have collected together a number of basic properties
of bounded packing in Section~\ref{sec:BoundedPacking},
and we state a variety of problems there for further thought.
In general, we would like to know:
\begin{ques}
Which countable groups have the property that all their
(finitely generated) subgroups have bounded packing?
\end{ques}

We also show in Section~\ref{sec:Cubulating} that under mild hypotheses
an action of a group $G$ on a $\CAT(0)$ cube complex $C$
gives rise to natural subgroups of $G$ with bounded packing.
These subgroups arise as stabilizers of hyperplanes in $C$.

\subsection{Pairwise close spaces}

Niblo--Reeves \cite{NibloReeves03} proved that given
$n$, $\kappa$, $D$ and $\delta$
there exists $M=M(n,\kappa,D,\delta)$
such that the following thinness condition holds:
If $A_1,\dots, A_n$ are $\kappa$--quasiconvex subspaces
of a $\delta$-hyperbolic space $X$
such that $d(A_i,A_j) < D$ for all $i,j$,
then there is a point $x\in X$ with $d(A_i,x)<M$ for all $i$.
This statement is also implicit in Sageev's work in \cite{Sageev97},
and plays a fundamental role in both these papers in proving the cocompactness
of the cube complex associated with a finite set of quasiconvex
codimension--$1$ subgroups.

We generalize this statement in Proposition~\ref{prop:CommonIntersection},
which contains a slightly stronger conclusion than
the following assertion that there is either a point
or a peripheral coset nearby.

\begin{prop}\label{prop:WeakCommonIntersection}
Let $G$ be a relatively hyperbolic group.
Choose positive constants $\kappa$ and $D$ and an integer $n \ge 1$.
Let $\Set{A}$ be an arbitrary set of $\kappa$--relatively quasiconvex subspaces of~$G$ such that
for each $A,A' \in \Set{A}$ the distance $d(A,A')$ is less than~$D$.
Then there is a constant $M=M(\kappa,D,n)$ so that at least
one of the following holds.
\begin{enumerate}
  \item For every set $\{A_1,\dots,A_n\}$ of $n$ distinct elements of $\Set{A}$,
  there is a point $x\in X$ such that
  $d(x,A_i) < M$ for each $i$.
  \item
  There is a peripheral coset $gP$ such that $d(gP,A) < M$
  for all $A \in \Set{A}$.
\end{enumerate}
\end{prop}

Proposition~\ref{prop:CommonIntersection} plays a significant role in our proof of
Theorem~\ref{thm:RelHypBoundedPacking}.

\subsection{Bounded Height and Width}
For relatively hyperbolic groups, the ``small'' subgroups
are either ``elliptic,'' meaning finite, or ``parabolic,''
which means conjugate into a peripheral subgroup.
In this context, we therefore redefine the \emph{height}
of $H$ in $G$ to be $(n-1)$ if
$n$ is the smallest number with the property that
for any $n$ distinct left cosets $g_1H,\dots, g_nH$,
the intersection $\bigcap_{1\leq i\leq n} g_i H g_i^{-1}$ is elliptic
or parabolic.
Similarly, the \emph{width} of $H$ in $G$ equals $(n-1)$
if $n$ is the smallest number with the property
that for any $n$ distinct left cosets $g_1H, \dots , g_nH$
the intersection $g_i H g_i^{-1} \cap g_j H g_j^{-1}$ is elliptic
or parabolic for some $i,j$.

The following is our main theorem about finite height and width:
\begin{thm}\label{thm:MainWidth}
Let $G$ be a relatively hyperbolic group, and let $H$ be a relatively quasiconvex subgroup.
Then $H$ has finite height.

Suppose, furthermore that $H \cap gPg^{-1}$ has bounded packing in $gPg^{-1}$
for each conjugate of each peripheral subgroup $P$. Then $H$ has finite width.
\end{thm}

We hope that further work in this area will resolve the following:
\begin{prob}
Let $G$ be relatively hyperbolic. Does every relatively quasiconvex subgroup of $G$
have finite width?
\end{prob}

If $H$ is a subgroup of $G$, the \emph{commensurator}
of $H$ in $G$ is the subgroup of all $g \in G$
such that $H \cap gHg^{-1}$ has finite index in both $H$
and $gHg^{-1}$.
We prove the following result on
commensurators of relatively quasiconvex subgroups
using ideas related to Theorem~\ref{thm:MainWidth}.

\begin{thm}
\label{thm:Commensurator}
Let $H$ be a relatively quasiconvex subgroup of a
relatively hyperbolic group $G$.
Suppose $H$ is infinite and nonparabolic.
Then $H$ has finite index in its commensurator.
\end{thm}

We use similar techniques to show that
a separable (relatively) quasiconvex subgroup $H$ of a (relatively)
hyperbolic group $G$ is (relatively) malnormal in a finite index subgroup $K$ of $G$ (see Theorem~\ref{thm:RelMalnormal}).
This result is new even in the hyperbolic case.

\subsection{Overview}

In Section~\ref{sec:BoundedPacking} we define the notion
of bounded packing and prove a number of basic results about
this concept.
We also give a number of examples and collect several problems.

In Section~\ref{sec:Cubulating} we discuss the relation between bounded packing
and actions on $\CAT(0)$ cube complexes.
By a result of Sageev, a codimension--$1$ subgroup $H < G$
with bounded packing gives rise to an action of $G$ on a finite dimensional $\CAT(0)$
cube complex.  In this section, we give a partial converse to Sageev's result.

In Section~\ref{sec:PackingHyperbolic} we prove the
bounded width theorem of Gitik--Mitra--Rips--Sageev \cite{GMRS98}
in the word hyperbolic setting.
The reader is advised
to understand this proof first, as our relatively hyperbolic generalization
will follow it closely, albeit with certain aspects substantially more complicated.

In Section~\ref{sec:RelHyp} we recount essential facts about
the geometry of the Cayley graph of a relatively hyperbolic group
with a finite generating set that were established
by \Drutu--Sapir in \cite{DrutuSapirTreeGraded}.
We use these results to prove two auxiliary results about triangles and quadrilaterals
with the property that each side lies near a peripheral subspace.

In Section~\ref{sec:RelativeCayley} we discuss the geometry of relative Cayley 
graphs, concentrating on the interplay between the usual (proper) metric and the
relative metric as applied to thinness conditions.

In Section~\ref{sec:RelQC} we analyze collections of relatively quasiconvex
subspaces in a relatively hyperbolic group. In particular,
we prove Proposition~\ref{prop:CommonIntersection}, an important technical 
result underlying this work that was described above.

In Section~\ref{sec:PackingRelHyp} we prove Theorem~\ref{thm:MainWidth},
our main result on bounded packing
and bounded width for a relatively quasiconvex subgroup of
a relatively hyperbolic group.
We also prove Theorem~\ref{thm:Commensurator}
on the commensurator of a relatively quasiconvex subgroup.

Section~\ref{sec:SeparableMalnormal} contains a brief proof of
the relative malnormality result
Theorem~\ref{thm:RelMalnormal} mentioned above,
using results from the previous section.

\subsection{Acknowledgements}

The first author's contribution to this material
is based upon work supported
by the National Science Foundation under Grant Nos.\ DMS-0505659,
DMS-0754254, and DMS-0808809.
The second author is supported by grants from NSERC.

\section{Bounded packing: First properties and examples}
\label{sec:BoundedPacking}

\begin{defn}[Bounded packing]
Let $G$ be a discrete group with a left invariant metric $d$.
Suppose also that $d$ is \emph{proper} in the sense that every
metric ball is finite. A subgroup $H$ has \emph{bounded packing}
in~$G$ (with respect to $d$) if, for each constant~$D$, there is a
number $N=N(G,H,D)$ so that for any collection of $N$ distinct
cosets $gH$ in~$G$, at least two are separated by a distance of at
least~$D$.
\end{defn}

A discrete group admitting a proper metric is necessarily countable.
On the other hand, if $G$ is finitely generated then
the word metric for any finite generating set is a proper, equivariant metric.
Restricting the word metric on $G$ to a subgroup $H$ gives a proper,
left invariant metric on $H$.  Since every countable group is a subgroup of a
finitely generated group, it follows that a discrete group $G$ admits a
proper, left invariant metric if and only if $G$ is countable.

\begin{lem}\label{lem:PackedIndependent}
If $G$ is countable, then bounded packing of a subgroup $H$ in $G$
is independent of the choice of proper, left invariant metric on $G$.
In other words, if $d_1,d_2$ are proper, left invariant metrics on $G$,
then $H$ has bounded packing in $G$ with respect to $d_1$
if and only if it has bounded packing with respect to $d_2$.
\end{lem}

\begin{proof}
Define $\rho \colon \N \to \N$ so that
$\rho(n)$ is the supremum of $d_1(1,g)$ over all $g$ with $d_2(1,g)<n$.
Note that there are only finitely many such $g$ since $d_2$ is proper,
so $\rho(n)$ is finite for each $n$.
By the left invariance of $d_1$ and $d_2$, it follows that
\[
   d_1(g_1,g_2) \le \rho\bigl( d_2(g_1,g_2) \bigr)
\]
for all $g_1,g_2 \in G$.

Now if $\Set{H}$ is a collection of left cosets of $H$
with pairwise $d_2$--distances at most $D$,
their pairwise $d_1$--distances are at most $\rho(D)$.
Thus if $H$ has bounded packing with respect to $d_1$,
it also has bounded packing with respect to $d_2$.
\end{proof}

The following result follows easily from the definition of bounded packing.

\begin{lem}\label{lem:PackedFiniteIndex}
Any finite index subgroup $K$ of a countable group $G$
has bounded packing in $G$.
\end{lem}

\begin{proof} The definition of bounded packing is vacuously satisfied
once $N$ is larger than $[G:H]$,
since there are no collections of $N$ distinct cosets of $H$ in $G$.
\end{proof}

\begin{lem}\label{lem:PackedPacked}
Suppose that $H\le K \le G$ and $G$ is countable.
\begin{enumerate}
\item \label{item:PackedInSubgroup}
If $H$ has bounded packing in $G$
then $H$ has bounded packing in $K$.
\item \label{item:PackedTransitive}
If $H$ has bounded packing in $K$, and $K$ has bounded packing in $G$
then $H$ has bounded packing in $G$.
\end{enumerate}
\end{lem}

\begin{proof}
We may assume the metric on $K$ is a restriction of the metric on $G$,
so that (\ref{item:PackedInSubgroup}) is immediate.

Let $\Set{H}$ is a collection of left cosets of $H$ in $G$ with pairwise
distances at most $D$.
The bounded packing of $H$ in $K$ bounds the number of elements of
$\Set{H}$ that lie in each left coset $gK$.
Furthermore, the bounded packing of $K$ in $G$ bounds the number of cosets
$gK$ that contain elements of $\Set{H}$.  Thus the size of $\Set{H}$
is bounded in terms of $D$, establishing (\ref{item:PackedTransitive}).
\end{proof}

\begin{prop}[Commensurability]\label{prop:CommensurablePacking}
Let $G$ be a countable group.
\begin{enumerate}
\item\label{item:FiniteIndexBottom}
Suppose $H \le K \le G$ and $[K:H]<\infty$.
Then $H$ has bounded packing in $G$ if and only if $K$ has bounded packing
in $G$.
\item\label{item:FiniteIndexTop}
Suppose $H \le K \le G$ and $[G:K]<\infty$.
Then $H$ has bounded packing in $K$ if and only if
$H$ has bounded packing in $G$.
\item\label{item:FiniteIndexIntersection}
Suppose $H,K \le G$ and $[G:K] < \infty$.
Then $H \cap K$ has bounded packing in $K$ if and only if $H$ has
bounded packing in $G$.
\end{enumerate}
\end{prop}

\begin{proof}
Let $d$ be a proper, left invariant metric on $G$.
To see (\ref{item:FiniteIndexBottom}) suppose $[K:H]<\infty$.
Then there is a constant $C$ such that in the metric $d$ each left coset $gK$
is at most a Hausdorff distance $C$ from a left coset $gH$
(the choice of $g$ is irrelevant).
Suppose $H$ has bounded packing in $G$. Let $\Set{K}$ be a collection
of left cosets of $K$ with pairwise distances at most $D$.
Replacing each coset $gK$ with a corresponding coset $gH$ gives a
collection $\Set{H}$ of left cosets of $H$ with the same cardinality as
$\Set{K}$ and with pairwise distances at most $D+2C$.
The bounded packing of $H$ in $G$ bounds the size of $\Set{H}$,
and hence also the size of $\Set{K}$, as a function of $D$.
Thus $K$ has bounded packing in $G$.

Conversely suppose $K$ has bounded packing in $G$.
By Lemma~\ref{lem:PackedFiniteIndex}, we know that $H$ has bounded packing
in $K$.  Therefore $H$ has bounded packing in $G$ by
Lemma~\ref{lem:PackedPacked}(\ref{item:PackedTransitive}),
completing the proof of (\ref{item:FiniteIndexBottom}).

Assertion (\ref{item:FiniteIndexTop}) follows immediately from
Lemmas \ref{lem:PackedFiniteIndex} and~\ref{lem:PackedPacked}.

Suppose now that $H,K \le G$ and $[G:K]\le \infty$.
Then $H \cap K$ has finite index in $H$.
Observe that (\ref{item:FiniteIndexIntersection}) is an immediate consequence
of (\ref{item:FiniteIndexBottom}) and~(\ref{item:FiniteIndexTop}),
since $H \cap K$ has bounded packing in $K$ if and only if $H \cap K$
has bounded packing in $G$ if and only if $H$ has bounded packing in $G$.
\end{proof}

\begin{cor}
\label{cor:FinitePacking}
If $G$ is countable, any finite subgroup $K \le G$ has bounded packing in $G$.
\end{cor}

\begin{proof}
Fix a proper, left invariant metric $d$ on $G$.
By Proposition~\ref{prop:CommensurablePacking}(\ref{item:FiniteIndexBottom})
it suffices to show that the trivial subgroup has bounded packing in $G$.
A left coset of the trivial subgroup is just a single element of $G$.
If $\Set{H}$ is any collection of pairwise $D$--close
left cosets of the trivial subgroup,
(in other words, elements of $G$) then
$\Set{H}$ is a subset of the ball of radius $D$
centered at any of its elements.
The number of elements in such a ball is finite since $d$ is proper,
and this number depends only on $D$ since $d$ is left-invariant.
\end{proof}

\begin{prop}\label{prop:PackingIntersection}
Suppose $H$ and $K$ have bounded packing in a countable group~$G$.
Then $H\cap K$ has bounded packing in $G$.
\end{prop}

\begin{proof}
Fix a proper, left invariant metric $d$ on $G$.
By Lemmas \ref{lem:PackedPacked}(\ref{item:PackedTransitive})
and~\ref{lem:PackedIndependent}
it suffices to show that $L:=H \cap K$ has bounded packing in $H$
with respect to $d$.
Let $\Set{L}$ be a set of left cosets of $L$ in $H$
whose pairwise $d$--distances are at most $D$.
If $hL$ and $h'L$ are distinct cosets of $L$ in $H$
then $hL = H \cap hK$ and $h'L = H \cap h'K$ for distinct cosets
$hK$ and $h'K$ in $G$.
Thus $\Set{L}$
is in one-to-one correspondence with a set $\hat{\Set{L}}$ of left cosets
of $K$ in $G$ whose pairwise distances are at most $D$.
But the bounded packing of $K$ in $G$ bounds the size of $\hat{\Set{L}}$
as a function of~$D$.
\end{proof}

\begin{lem}\label{lem:PackingPreimage}
Let $1\rightarrow N\rightarrow G\rightarrow Q\rightarrow 1$ be a short exact
sequence of countable groups.
Let $H$ be a subgroup of $G$ which projects to the
subgroup $\bar H$ of $Q$.
Then $\bar H$ has bounded packing in $Q$
if and only if $HN$ has bounded packing in $G$.
\end{lem}

\begin{proof}
For each $g \in G$, let $\bar g$ denote the image of $g$ in $Q$.
Fix a proper, left invariant metric $d$ on $G$.
Define $\bar d \colon Q \to \R$ so that $\bar{d}(q_1,q_2)$ is the infimum of
$d(g_1,g_2)$ over the set of all $g_1,g_2$ with $q_i = \bar g_i$.
A straightforward argument shows that $\bar{d}$ is a proper,
left invariant metric on $Q$.

The projection $G \to Q$ induces a one-to-one correspondence between
left cosets of $HN$ in $G$ and left cosets of $\bar H$ in $Q$.
It is clear that $d(g_1,g_2) \ge \bar{d}(\bar{g_1},\bar{g_2})$
for all $g_1,g_2 \in G$.
Therefore
\[
   d(xHN,yHN) \ge \bar{d}(\bar x \bar H, \bar y \bar H)
\]
for all $x,y \in G$.

Conversely, any element $x \in Q$ with $\bar{d}(x,\bar{1})=a$
lifts to an element $y \in G$ with $d(y,1)=a$.
If $\bar{d}(\bar{g_1}\bar H, \bar{g_2}\bar H) = n$
then there is an element $x \in Q$ with $\bar{d}(x,\bar{1}) = a$
such that $\bar{g_1}\bar{H}\bar{w}$ intersects $\bar{g_2}\bar{H}$.
If $y$ is a lift of $x$ with $d(y,1)=a$,
then $ g_1 H y N = g_1 HN y$ intersects $ g_2 HN$,
so $d(g_1 HN, g_2 HN) \le a$.

Therefore the distance between left cosets of $HN$ in $G$ is equal to the
distance between corresponding left cosets of $\bar H$ in $Q$.
Evidently $HN$ has bounded packing in $G$ if and only if
$\bar H$ has bounded packing in $Q$.
\end{proof}

\begin{cor}\label{cor:PackingNormal}
Each normal subgroup $N$
of a countable group $G$ has bounded packing.
\end{cor}

\begin{proof}
Apply Lemma~\ref{lem:PackingPreimage} with $H$ equal to the trivial subgroup.
\end{proof}

\begin{rem}
In fact, a normal subgroup $N$ has the stronger property that
for each $D$ there exists $M=M(D)$ such that
at most $M$ left cosets $gN$ satisfy $d(N,gN)< D$.
\end{rem}

\begin{cor}
Every subgroup of a countable abelian group has bounded packing.\qed
\end{cor}

\begin{thm}
Let $N$ be a countable, virtually nilpotent group.
Then each subgroup of $N$ has bounded packing in $N$.
\end{thm}
We note that this can be also be proven using
either Lemma~\ref{lem:central extension} or Corollary~\ref{cor:subnormal}.

\begin{proof}
By
Proposition~\ref{prop:CommensurablePacking}(\ref{item:FiniteIndexIntersection})
it suffices to prove the theorem
when $N$ is a countable nilpotent group.
The theorem is obvious if $N$ is the trivial group.
We proceed by induction on the length of the lower central series.
Let $H$ be a subgroup of $N$.
Let $Z$ denote the center of $N$.
Let $\bar N := N / Z$ and let $\bar H$ be the image of $H$ in $\bar N$.
By induction,
$\bar H$ has bounded packing in $\bar N$.
Therefore by Lemma~\ref{lem:PackingPreimage},
$HZ$ has bounded packing in $N$.
Because $Z$ is central, $H$ is normal in $HZ$
and so $H$ has bounded packing in $HZ$ by Corollary~\ref{cor:PackingNormal}.
Finally, by Lemma~\ref{lem:PackedPacked}(\ref{item:PackedTransitive}),
$H$ has bounded packing in $G$ since,
$H$ has bounded packing in $HZ$ and $HZ$ has bounded packing in $G$.
\end{proof}

\begin{ques}\label{ques:solvable}
Let $S$ be a solvable group.
Does every subgroup of $S$ have bounded packing in $S$?
\end{ques}

It seems too much to expect an affirmative answer to Question~\ref{ques:solvable},
however we do expect the following to hold:
\begin{conj}
Let $P$ be virtually polycyclic.
Then each subgroup of $P$ has bounded packing in $P$.
\end{conj}

The following is an immediate consequence of Corollary~\ref{cor:PackingNormal}
together with Lemma~\ref{lem:PackedPacked}(\ref{item:PackedTransitive}).

\begin{cor}\label{cor:subnormal}
Let $H$ be a subgroup of a countable group $G$.
If $H$ is subnormal then $H$ has bounded packing in $G$.\qed
\end{cor}

\begin{exmp}
For $n\neq 0$, let $BS(1,n)$ denote the group presented by $\presentation{a,t}{tat^{-1}=a^n}$.
It is well-known that for $n\neq \pm1$, the subgroup $\langle a \rangle$
is not quasi-isometrically embedded;
in fact, it is exponentially distorted.
Nevertheless $\langle a \rangle$ has bounded packing in $BS(1,n)$.
Indeed, there is a short exact sequence
$1\rightarrow A \rightarrow BS(1,n) \rightarrow \Z \rightarrow 1$
with $\gen{a} \le A$ such that $A$ is isomorphic to $\Z[1/n]$.
Since $A$ is abelian, $\gen{a}$ is subnormal in $BS(1,n)$ and hence has bounded
packing by Corollary~\ref{cor:subnormal}.

It is interesting to note that $A$ is not finitely generated.  (Indeed, $A$ is 
locally cyclic but not itself cyclic.)
Thus consideration of bounded packing in non--finitely generated groups
can give us useful information about finitely generated subgroups of a finitely generated group.
\end{exmp}


As the following Lemma shows, the group $F_2\times \integers$
has bounded packing with respect to each of its finitely generated subgroups.

\begin{lem}\label{lem:central extension}
Let $1\rightarrow C \rightarrow G \rightarrow Q \rightarrow 1$
be a short exact sequence with $C$ central in the countable group $G$.
If every \textup{[}finitely generated\textup{]} subgroup of $Q$ has bounded packing
then every \textup{[}finitely generated\textup{]} subgroup of $G$ has bounded packing.
\end{lem}
\begin{proof}
Let $H$ be a [finitely generated]
subgroup of $G$. By hypothesis, its image $\bar H$ has bounded packing in $Q$.
Now $H$ has bounded packing in $CH$ by Corollary~\ref{cor:PackingNormal},
and $CH$ has bounded packing in $G$ by Lemma~\ref{lem:PackingPreimage}.
Therefore $H$ has bounded packing in $G$ by Lemma~\ref{lem:PackedPacked}.
\end{proof}

\begin{exmp}
Let $M$ be a Seifert-fibered $3$--manifold.
Then there is a short exact sequence
$1\rightarrow \integers \rightarrow \pi_1M \rightarrow \pi_1S\rightarrow 1$,
where $S$ is a $2$--dimensional orbihedron.
Since every finitely generated subgroup of $\pi_1S$ has bounded packing
(by local quasiconvexity and Theorem~\ref{thm:HyperbolicPacking} below)
it follows from Lemma~\ref{lem:central extension}
that every finitely generated subgroup of $\pi_1M$ has bounded packing.
\end{exmp}

Let us now examine subgroups that do not have bounded packing.
At present suprisingly few such examples are known.
In fact, the authors know of only one basic example, due to Rubinstein--Wang.
The group is the fundamental group of a $3$--dimensional graph manifold;
that is, a manifold formed by gluing Seifert fibered manifolds together
along tori.
In light of the preceding example,
we see that the bounded packing property is not preserved
by very simple graphs of groups.

\begin{exmp}
\label{exmp:RubinsteinWang}
Rubinstein--Wang \cite{RubinsteinWang98} constructed a finitely generated
surface subgroup of the fundamental group of a graph manifold
that doesn't have bounded packing.
Their example is an immersed $\pi_1$--injective
surface $S$ in a $3$--dimensional graph manifold $M$.
Let $\widetilde S\rightarrow \widetilde M$ denote a lift
of $S\rightarrow M$ to a map between the universal covers.
The various lifts correspond to translates $g\widetilde S$
where $g$ varies over the various cosets $g\pi_1S$.
Remarkably, any two such lifts $g_i\widetilde S$ and $g_j\widetilde S$
intersect each other.
Consequently any two cosets $g_i\pi_1S$ and $g_j\pi_1S$ are a uniformly
bounded distance from each other in $\pi_1M$.
\end{exmp}

It seems reasonable to expect subgroups without bounded packing
to appear in the fundamental groups of other graph manifolds.
The following problem is a reasonable test case for this phenomenon.

\begin{prob}
Determine exactly which subgroups of the following group $G$ have bounded
packing:
\[
   G:=\presentation{a,b,c,d}{[a,b], [b,c], [c,d]}
\]
(The group $G$ is the fundamental group of a graph manifold.)
\end{prob}

Since there are subgroups of countable groups without bounded packing,
it follows from Lemma~\ref{lem:PackingPreimage}
that (non--finitely generated) subgroups of the free group $F_2$
do not always have bounded packing.
Notice that,
since finitely generated subgroups of $F_2$ are quasiconvex,
Theorem~\ref{thm:HyperbolicPacking} below implies that
finitely generated subgroups of $F_2$ always have bounded packing.

Rips gave the following construction in \cite{Rips82}:
\begin{prop}\label{prop:rips}
Let $Q$ be a finitely presented group.
Then there exists a short exact sequence
$1\rightarrow N \rightarrow G \rightarrow Q \rightarrow 1$
such that $N$ is finitely generated,
and $G$ is a finitely presented $C'\!\bigl(\frac{1}{6}\bigr)$ group.
\end{prop}

\begin{exmp}\label{exmp:rips short exact}
Using Rips's construction, we will construct a word hyperbolic group $G$
and a finitely generated subgroup $K$ such that $K$ does not have
bounded packing in $G$.
Let $Q$ be a finitely presented group with a subgroup $H$
that fails to have bounded packing in $Q$.
Let $1\rightarrow N \rightarrow G \rightarrow Q \rightarrow 1$
be the short exact sequence from Proposition~\ref{prop:rips}.
Let $K$ be the preimage of $H$ in $G$, and note that
$K$ is finitely generated since both $N$ and $H$ are.
By Lemma~\ref{lem:PackingPreimage},
$K$ does not have bounded packing in $G$.
Observe that $G$ is word hyperbolic, since it is finitely presented and
$C'\!\bigl(\frac{1}{6}\bigr)$.
\end{exmp}

\begin{rem}
We can choose $H$ to be a codimension--$1$ subgroup without bounded
packing to make the group $G$ in Example~\ref{exmp:rips short exact}
have an exotic action on an infinite dimensional cube complex.
\end{rem}

\begin{prob}
Give an example of a cyclic subgroup $Z$ of a
finitely generated group $G$
such that $Z$ does not have bounded packing in $G$.
\end{prob}

\begin{ques}
Are there finitely generated subgroups of $F_2\times F_2$
without bounded packing?
\end{ques}

\begin{ques}
Let $1 \rightarrow N \rightarrow G \rightarrow A \rightarrow 1$ be
an extension of finitely generated groups.
Suppose every subgroup of $N$ has bounded packing in $N$, and $A$ is abelian.
Does every subgroup of $G$ have bounded packing in $G$?
\end{ques}

\begin{prob}
Let $H$ be a retract of the finitely generated group $G$.
Does $H$ have bounded packing in $G$?
\end{prob}

\section{Bounded packing and actions on $\CAT(0)$ cube complexes}
\label{sec:Cubulating}

Suppose $G$ is a group with a finite generating set $\Set{S}$.
A subgroup $H < G$ is \emph{codimension}--$1$ if there is a constant $L>0$
so that the $L$--neighborhood of $H$ in $\Cayley(G,\Set{S})$
has at least two ``deep'' complementary components.
A complementary component is \emph{deep} if it contains elements whose distance from
$H$ is arbitrarily large.

A theorem of Sageev \cite{Sageev95}, together with a result proved independently by Gerasimov \cite{Gerasimov97} and
Niblo--Roller \cite{NibloRoller98},
shows that a group $G$ with a codimension--$1$ subgroup $H$
has an action on a $\CAT(0)$ cube complex $C$ without a global fixed point.
All hyperplanes in $C$ lie in a single orbit under the action of $G$,
and $H$ is the stabilizer of some hyperplane $\Lambda$.
The dimension of $C$ is the maximal size of a set of pairwise
transverse hyperplanes.
The relation to bounded packing is the following basic fact, observed by Sageev \cite{Sageev97}:
If $H$ is finitely generated and
hyperplanes $g_0 \Lambda$ and $g_1 \Lambda$ in $C$
are transverse, then $d_{\Set{S}} (g_0 H, g_1 H) < M$
for some universal constant $M$.

\begin{cor}[Sageev]
Suppose $H$ is a finitely generated
codimension--$1$ subgroup of a finitely generated group $G$.
If $H$ has bounded packing in $G$, then the corresponding $\CAT(0)$ cube 
complex $C$ is finite dimensional.
\end{cor}

The goal of this section is to prove the following converse, of sorts, to the preceding result.

\begin{thm}
\label{thm:CubicalPacking}
Suppose the countable discrete group $G$ acts on a $\CAT(0)$ cube complex $C$,
and $H\le G$ is the stabilizer of a hyperplane $\Lambda$.
\begin{enumerate}
\item \label{item:LocFinite}
If $C$ is locally finite, then $H$ has bounded packing in $G$.
\item \label{item:CocompactStab}
Let $V(\Lambda)$ be the set of vertices incident to the edges that
cross $\Lambda$.
If $G$ acts metrically properly on $C$,
and $H$ acts on $V(\Lambda)$ with only finitely many orbits of vertices,
then $H$ has bounded packing in $G$.
\end{enumerate}
\end{thm}

In a sense this theorem is optimal, as shown by the following example.

\begin{exmp}
Recall that Example~\ref{exmp:RubinsteinWang}
gives a group $G$ and a subgroup $H$ such that $H$ does
not have bounded packing in $G$.  The group $H$ is the fundamental group
of an immersed surface $S$
in a $3$--dimensional graph manifold $M$ with $G = \pi_1(M)$.
Since $H$ is a codimension--$1$ subgroup of $G$,
it follows that $G$ has an action on a $\CAT(0)$ cube complex $C$ such that
$H$ is the stabilizer of a hyperplane $\Lambda$.
Condition~(\ref{item:LocFinite}) fails resolutely since $C$
is an infinite dimensional cube, as any two translates of the surface
$S$ cross.
The authors believe failure of (\ref{item:CocompactStab})
can be traced to a failure of the cocompactness of the action
of $H$ on $V(\Lambda)$.
Indeed even if the action of $G$ on $C$ is not proper,
the authors
believe this can be remedied with the addition of further codimension--$1$
subgroups leading to a metrically proper action on a new cube complex $C'$.
The action of $H$ on the new hyperplane $\Lambda'$ would necessarily
fail to be cocompact by Theorem~\ref{thm:CubicalPacking},
as $H$ does not have bounded packing in $G$.
\end{exmp}

Theorem~\ref{thm:CubicalPacking}
follows from several results about convex subsets of a
$\CAT(0)$ cube complex $C$.
Let $C^{1}$ denote the $1$--skeleton of a $\CAT(0)$ cube complex with its induced path metric $d$.
It is well-known that the metric $d$ on the vertices coincides with the ``wall metric.''
That is to say, $d(v,w)$ is equal to the number of hyperplanes separating $v$ from $w$.
The \emph{interval} $[v,w]$ is the set of all vertices $u$ that lie on some $d$--geodesic from
$v$ to $w$.
A subset $S \subseteq C^{0}$ is \emph{$d$--convex} if every interval connecting
two elements of $S$ is contained in $S$.

We also find it useful to consider the graph $C^{\Delta}$ obtained from $C^{1}$
by including an edge between two vertices if they lie in a common cube.
Let $d^{\Delta}$ be the path metric on $C^{\Delta}$.

\begin{lem}
\label{lem:Cubical1Neighborhood}
Let $C$ be a $\CAT(0)$ cube complex, and let $[r,s]$ and $[t,u]$ be intervals
such that $t$ and $u$ lie in the $1$--neighborhood of $[r,s]$ in $C^\Delta$.
Then each vertex of $[t,u]$ lies in the $1$--neighborhood of $[r,s]$ in $C^\Delta$
as well.
\end{lem}

\begin{proof}
It is enough to show that any geodesic $c$ from $t$ to $u$ lies in the
$1$--neighborhood of $[r,s]$ in $C^\Delta$
if its endpoints lie in the $1$--neighborhood of $[r,s]$ in $C^\Delta$.
If $d(t,u) \le 1$, the result is trivial, so we induct on
$\ell:=d(t,u)$ for $\ell\ge 2$.

By hypothesis, there is a cube $A$ containing $t$ and intersecting $[r,s]$.
Let us choose $A$ minimal with this property.
Then the hyperplanes $\Lambda_1,\dots,\Lambda_n$ separating $t$ from $[r,s]$
are the same as the hyperplanes transverse to $A$.
(In particular, $A$ is an $n$--cube.)

Let $v$ be the vertex of $c$ adjacent to $t$, and let $\Pi$ be the hyperplane separating
$v$ from $t$.
If $\Pi=\Lambda_i$ for some $i$, then $v$ is a vertex of the cube $A$,
and we are done by induction.

If $\Pi$ is transverse with $\Lambda_i$ for all $i$,
then the collection $\{\Pi,\Lambda_1,\dots,\Lambda_n\}$ is pairwise transverse.
Since $\Pi$ and the $\Lambda_i$ are each transverse to edges incident to $t$, there is an
$(n+1)$--cube $A'$ containing $A$ as a subcube,
such that $\Pi$ and each $\Lambda_i$ are transverse to $A'$.
Since $A$ intersects $[r,s]$ and $A'$ contains $v$,
the result follows from the inductive hypothesis.

The only remaining possibility is that $\Pi$ is nested with some $\Lambda_i\ne\Pi$.
Then $\Lambda_i$ separates $[r,s]$ from $\{t,v\}$
and $\Pi$ separates $[r,s] \cup \{t\}$ from $\{v\}$.
In this case, any path in $C^{1}$ from $v$ to $[r,s]$ must cross
$\Pi$ first and then $\Lambda_i$.
Since $c$ is a geodesic in the wall metric, it crosses $\Pi$ only once,
so $\Pi$ lies in the pairwise transverse family $\{\Pi_1,\dots,\Pi_m\}$ of hyperplanes
separating $u$ from $[r,s]$, corresponding to a minimal cube $B$
that contains $u$ and intersects $[r,s]$.
But then $\Lambda_i$ is also among the $\Pi_j$, contradicting our assumption that
$\Lambda_i$ and $\Pi$ are nested.
\end{proof}

\begin{cor}\label{cor:ConvexNbd}
Let $C$ be a $\CAT(0)$ cube complex, and let $K$ be a $d$--convex subset of $C^{0}$.
For each $n\ge 0$, the $n$--neighborhood of $K$ in $C^\Delta$
is also $d$--convex.
\end{cor}

\begin{proof}
It is clearly sufficient to prove the result when $n=1$, since the general result follows by
induction on $n$.
But the case $n=1$ is an immediate consequence of Lemma~\ref{lem:Cubical1Neighborhood}.
\end{proof}

A collection of subsets $\Set{F}$ of a space has the \emph{Helly property}
if, whenever $\Set{F}_0\subseteq\Set{F}$ is a finite collection of pairwise intersecting
subsets, the total intersection $\bigcap_{F\in\Set{F}_0} F$ is nonempty.

\begin{cor}
If $C$ is any $\CAT(0)$ cube complex, then the family of $d$--convex subsets
of $C^{0}$ has the Helly property.
\end{cor}

\begin{proof}
It has been shown by Gerasimov \cite{Gerasimov97},
Roller \cite{RollerPocSets},
and Chepoi \cite{Chepoi00}, independently, that the $1$--skeleton of a
$\CAT(0)$ cube complex is a median graph,
meaning that for any three vertices $x$, $y$ and $z$,
the intervals $[x,y]$, $[y,z]$ and $[x,z]$
have triple intersection consisting of a single vertex.
In a median graph, it is clear from the definition
that any collection of three pairwise intersecting convex sets has a nonempty triple intersection.
The Helly property now follows by an elementary induction argument.
\end{proof}

The previous result generalizes a result of Bandelt and
van~de~Vel \cite{BandeltVanDeVel91}
in the setting of
median graphs stating that the collection of all balls in $C^\Delta$ has the 
Helly property.
Note that $d^\Delta$--balls are $d$--convex by Corollary~\ref{cor:ConvexNbd}.
Graphs in which the family of balls has the Helly property are known as
\emph{Helly graphs}.

We also record the following well-known fact about hyperplanes,
which we derive from results above.

\begin{prop}
Let $\Lambda$ be a hyperplane in a $\CAT(0)$ cube complex,
and let $V(\Lambda)$ be the set of vertices of the edges crossing $\Lambda$.
Then $V(\Lambda)$ is $d$--convex.
\end{prop}

\begin{proof}
Each halfspace of $\Lambda$ is clearly $d$--convex.  But $V(\Lambda)$
is the set of points within a $d^\Delta$--distance $1$ of both halfspaces of $\Lambda$.
Since $d^\Delta$--neighborhoods of $d$--convex sets are convex
and intersections of $d$--convex sets are $d$--convex, the result is obvious.
\end{proof}

We are now ready to prove Theorem~\ref{thm:CubicalPacking}.

\begin{proof}[Proof of Theorem~\ref{thm:CubicalPacking}]
Suppose $H$ stabilizes the hyperplane $\Lambda$ in $C$.
Let $V(\Lambda)$ be the convex set of vertices of the edges that cross $\Lambda$.
Fix a proper, left invariant metric $d_G$ on $G$ and a constant $D>0$.
Suppose $\Set{H}$ is a finite collection of cosets of $H$
whose pairwise distances in $G$ are at most $D$.

Choose a basepoint $x \in V(\Lambda)$, and let
\[
   D' := \sup d^\Delta\bigl(x,g(x)\bigr)
\]
where the supremum ranges over the finite set of group elements whose $d_G$--distance from $1$
is at most $D$.
If cosets $gH$ and $g'H$ are at a distance at most $D$ in $G$, then the sets
$gH(x)\subset V(g\Lambda)$ and $g'H(x)\subset V(g'\Lambda)$ are at a $d^\Delta$--distance
at most $D'$.  Therefore the translates $V(g\Lambda)$ and $V(g'\Lambda)$ are at a
$d^\Delta$--distance at most $D'$ as well.  It follows immediately that the
$D'$--neighborhoods of the translates $V(g\Lambda)$ are pairwise intersecting $d$--convex
sets in $C$.  Hence there is a point $p$ within a $d^\Delta$--distance $D'$
of every $V(g\Lambda)$.

If $C$ is uniformly locally finite, then there is a uniform bound on the number of
hyperplanes intersecting any $d^\Delta$--ball of radius $D'$ in $C$.
Hence the cardinality of $\Set{H}$ is bounded as well.

Now suppose $G$ acts metrically properly on $C$, and $H$ acts on $V(\Lambda)$ with a finite quotient.
Then $V(\Lambda)$ lies in the $R$--neighborhood of the orbit $Hx$ in $C^\Delta$ for some $R>0$.
Thus the $d^\Delta$--ball centered at $p$ with radius $D'+R$ intersects $gH(x)$
for each coset $gH \in \Set{H}$.
Since $G$ acts metrically properly on $C$, there is a finite upper bound on the size of $\Set{H}$,
as desired.
\end{proof}

\section{Bounded packing in hyperbolic groups}
\label{sec:PackingHyperbolic}

In this section, we give a self-contained proof in the word hyperbolic
setting that quasiconvex subgroups have bounded packing.

\begin{defn}
Let $H\leq G$ be a subgroup. The \emph{height} of $H$ in $G$,
denoted $\height_G(H)$, is the maximal number of distinct cosets
$g_iH$ so that the intersection $\bigcap_i g_i H g_i^{-1}$ is
infinite. The \emph{width} of $H$ in $G$, denoted $\width_G(H)$,
is the maximal number of distinct cosets $g_iH$ so that for all
$i,j$ the intersection $g_i H g_i^{-1} \cap g_j H g_j^{-1}$ is
infinite. In case there is no maximum then we say that the
height/width is infinite.
\end{defn}

Finally, note that
the height and width of a finite subgroup are both equal to zero.

\begin{lem}\label{lem:HeightDecrease}
Suppose $H \le G$ has height $0< n < \infty$ in~$G$.
Choose $g \in G$ so that $gH \ne H$, and let $K := H \cap gHg^{-1}$.
Then $\height_H(K) < n$.
\end{lem}

\begin{proof}
Choose $h_1, \dots, h_n \in H$ so that the cosets $h_iK$ are distinct.
We will show that the intersection $\bigcap h_i K h_i^{-1}$ is finite.
Note that
\[
   \bigcap h_i K h_i^{-1} = \bigcap h_i (H \cap gHg^{-1}) h_i^{-1}
        = H \cap \left( \bigcap h_ig H g^{-1}h_i^{-1} \right)
\]
is an intersection of $n+1$ conjugates of~$H$ in~$G$.
Since $\height_G(H) =n$, it suffices to show that
the elements $1,h_1g, \dots, h_ng$ represent distinct left cosets of~$H$
in~$G$.

First note that $H \ne h_igH$ since $H \ne gH$.
Now suppose two cosets $h_igH$ and $h_jg H$ are equal.
Then we have $gH = h_i^{-1}h_j gH$, which implies that
$h_i^{-1}h_j \in gHg^{-1}$.  But $h_i^{-1}h_j \in H$ as well,
so $h_i^{-1}h_j \in H \cap gHg^{-1} = K$.
Thus $h_i K = h_j K$, and we must have $i=j$.
\end{proof}

\begin{lem}\label{lem:LineNearQC}
Let $Y$ be a $\kappa$--quasiconvex subspace of a $\delta$--hyperbolic
space~$X$.  Then there exists $\kappa' = \kappa'(\kappa,\delta)$
so that any geodesic line~$A$ in a finite neighborhood of~$Y$
lies in a $\kappa'$--neighborhood of~$Y$.
\end{lem}

\begin{proof}
Let $c$ be an arbitrary line in a finite neighborhood of~$Y$.
Then the endpoints $c^+$ and $c^-$ of $c$ lie in the limit set of~$Y$.
Let $c'$ be a line connecting $c^+$ and $c^-$ that is
a limit of geodesic segments connecting points of~$Y$.
In a $\delta$--hyperbolic space, two geodesic lines with the same
endpoints at infinity are at a Hausdorff distance at most $2\delta$.
But by quasiconvexity, $c' \subseteq \nbd{Y}{\kappa}$,
so $c \subseteq \nbd{Y}{\kappa+2\delta}$.
\end{proof}

\begin{lem}\label{lem:FiniteHeight}
Let $H$ be a $\kappa$--quasiconvex subgroup of a $\delta$--hyperbolic
group~$G$.  Then $\height_G(H)$ is finite.
\end{lem}

\begin{proof}
Let $g_1H, \dots g_n H$ be distinct cosets in~$G$, and suppose
$\bigcap g_i H g_i^{-1}$ contains an infinite order (hyperbolic)
element~$x$.
Let $A$ be a geodesic axis for~$x$.  Then $A$ lies in a finite
neighborhood of each coset $g_iH$, so
by Lemma~\ref{lem:LineNearQC} the size of these neighborhoods is uniformly
bounded by~$\kappa'$.
In particular, each coset $g_i H$ intersects the ball of radius~$\kappa'$
about some point $a \in A$.  But there is a uniform bound~$N$
on the number of cosets of~$H$ intersecting any metric ball
of radius~$\kappa'$.
\end{proof}

\begin{lem}\label{lem:CloseCosets}
Let $G$ have a finite generating set~$\Set{A}$, and suppose $xH$ and $yK$
are arbitrary left cosets of subgroups of~$G$.
For each constant $L$ there is a constant $L'=L'(G,\Set{A},xH,yK)$
so that in $G$ with the word metric $d_\Set{A}$ we have:
\[
   \nbd{x H}{L} \cap \nbd{y K}{L} \subseteq
     \nbd{x H x^{-1} \cap y K y^{-1}}{L'}
\]
\end{lem}

\begin{proof}
%
%
%
If there is no such $L'$, then there is a sequence $(z_i)$ in~$G$ so that
$z_i$ is in the $L$--neighborhood of both $xH$ and $yK$,
but $i < d(z_i, xHx^{-1} \cap yKy^{-1})$ for each~$i$.
It follows that $z_i = x h_i p_i = y k_i q_i$
for some $h_i\in H$, $k_i \in K$ and $p_i,q_i \in G$
with $\abs{p_i}_{\Set{A}}$ and $\abs{q_i}_{\Set{A}}$ both less than~$L$.
Passing to a subsequence if necessary, we may assume that $p_i$
and $q_i$ are constants $p$ and $q$, so that
for each $i$ we have $z_i = x h_i p = y k_i q$.
Therefore
\[
   z_i z_1^{-1} = x h_i h_1^{-1} x^{-1} = y k_i k_1^{-1} y^{-1}
     \in xHx^{-1} \cap yKy^{-1}.
\]
It follows that the distance between $z_i$
and $xHx^{-1} \cap yKy^{-1}$ is at most
$\abs{z_1}_{\Set{A}}$ for all~$i$, contradicting our choice of $(x_i)$.
\end{proof}

The following result follows immediately by induction from
Lemma~\ref{lem:CloseCosets}.

\begin{cor}
\label{cor:ManyCloseCosets}
Let $G$ be a finitely generated group with a fixed generating set $\Set{A}$.
For each finite set of cosets $g_1 H_1, \dots g_n H_n$ and each $M>0$
there is a constant $M'>0$ such that:
\[
   \nbd{g_1 H_1}{M} \cap \dots \cap \nbd{g_n H_n}{M}
     \subseteq \nbd{ g_1 H_1 g_1^{-1}
              \cap \dots \cap g_n H_n g_n^{-1} }{M'}
    \rlap{\hspace{.5in}\qedsymbol}
\]
\end{cor}

\begin{lem}\label{lem:QCNeighborhood}
Let $X$ be a $\delta$--hyperbolic space, and let $Y$ be a
$\kappa$--quasiconvex subspace.
For any $\xi>0$ there is an $\eta = \eta(\delta,\kappa,\xi)$ so that the
$\xi$--neighborhood of~$Y$ is $\eta$--quasiconvex.
\end{lem}

\begin{proof}
Choose points $x,y \in \nbd{Y}{\xi}$, and $z,w \in Y$
so that $d(x,z)$ and $d(y,w)$ are less than~$\xi$.
Since quadrilaterals in $X$ are $2\delta$--thin,
the geodesic $[x,y]$ lies in the $2\delta$--neighborhood
of $[x,z] \cup [z,w] \cup [w,y]$.
But this piecewise geodesic lies in the
$(\kappa+\eta)$--neighborhood of~$Y$.
\end{proof}

\begin{thm}[cf Gitik--Mitra--Rips--Sageev \cite{GMRS98}]
\label{thm:HyperbolicPacking}
Let $H$ be a quasiconvex subgroup of a $\delta$--hyperbolic
group~$G$.
Then $H$ has bounded packing in $G$.
\end{thm}

Before giving the proof of Theorem~\ref{thm:HyperbolicPacking},
we record the following:

\begin{cor}[Gitik--Mitra--Rips--Sageev \cite{GMRS98}]
If $H$ is a quasiconvex subgroup of a $\delta$-hyperbolic group~$G$,
then $\width_G(H)$ is finite.
\end{cor}

\begin{proof}
As in Lemma~\ref{lem:FiniteHeight}, apply Lemma~\ref{lem:LineNearQC}
to see that any two conjugates with infinite intersection
have cosets uniformly close together.
\end{proof}

\begin{proof}[Proof of Theorem~\ref{thm:HyperbolicPacking}]
Fix finite generating sets $\Set{A}$ and $\Set{B}$ for $G$ and $H$
respectively.
By Lemma~\ref{lem:FiniteHeight}, we know that $\height_G(H)$ is finite.
We will prove the theorem by induction on height.
The case where $H$ has height zero is easy, since in this case $H$ is
finite and the result follows from Corollary~\ref{cor:FinitePacking}.

Now assume by induction that the theorem holds for every
hyperbolic group $G'$ and quasiconvex subgroup $H'$
with $\height_{G'}(H') < \height_G(H)$.
Let $\Set{H}$ be a set of left cosets $gH$ whose pairwise
distances are at most~$D$.  Our goal is to show that $\Set{H}$ is
finite and to bound the cardinality of~$\Set{H}$ as a function of~$D$.
Translating $\Set{H}$ if necessary, we may assume that $H \in \Set{H}$.
Observe that if $d(gH,H) < D$ then the coset $hgH$
is within a distance $D$ of the identity for some $h \in H$.
Since the metric $d_{\Set{A}}$ is proper on $G$,
the ball of radius $D$ centered at the identity is finite.
It follows that the left cosets $gH$ intersecting $\nbd{H}{D}$ lie in at most
$f$ distinct $H$--orbits for some $f=f(D)<\infty$.
Thus it suffices to bound the number of elements of $\Set{H}$
in the orbit $H(gH)$ for each fixed $g \notin H$.

If we let $K := H \cap gHg^{-1}$,
then Lemma~\ref{lem:HeightDecrease} shows that
\[
   \height_H(K) < \height_G(H).
\]
Since $K$ is a quasiconvex subgroup of the hyperbolic group $H$,
the inductive hypothesis applied
to $K\le H$ gives for each $D'$ a number $M' = M'(D') < \infty$
so that any collection of $M'$ distinct cosets $hK$ in $H$ contains a pair
separated by a $\Set{B}$--distance at least~$D'$.
Furthermore, the proof of Lemma~\ref{lem:HeightDecrease}
shows that there is a well-defined map $hgH \to hK$ taking
left cosets of $H$ in the orbit of $gH$ to left cosets of~$K$.
A similar argument shows that this map is bijective.

In order to complete the proof, we will show that $D$--closeness of distinct
cosets $h_0gH$ and $h_1gH$ in $(G,d_{\Set{A}})$ implies $D'$--closeness of
the corresponding cosets $h_0K$ and $h_1K$ in $(H,d_{\Set{B}})$
for some $D'$ depending on~$D$.
Roughly speaking, this claim is proved by considering a $\delta$--thin
triangle in~$G$ whose sides are close to the three cosets
$H$, $h_0gH$ and $h_1gH$, and applying Lemma~\ref{lem:CloseCosets}.
More precisely, suppose we have points $x$, $y$ and~$z$
such that
\begin{align*}
   x &\in \nbd{H}{D} \cap \nbd{h_0gH}{D}, \\
   y &\in \nbd{H}{D} \cap \nbd{h_1gH}{D} \quad \text{and} \\
   z &\in \nbd{h_0gH}{D} \cap \nbd{h_1gH}{D}.
\end{align*}
By Lemma~\ref{lem:QCNeighborhood} we have
\begin{align*}
   [x,y] \subset \nbd{H}{D_1}, \quad [x,z] \subset \nbd{h_0gH}{D_1}
   \quad \text{and} \quad
   [y,z] \subset \nbd{h_1gH}{D_1}
\end{align*}
for some $D_1$ depending on~$D$.
By $\delta$--hyperbolicity of $(G,d_{\Set{A}})$,
there is a point $w \in [x,y]$
within a distance $\delta$ of both $[x,z]$ and $[y,z]$.
Applying the isometry $h_i$ to the conclusion of Lemma~\ref{lem:CloseCosets},
we see that
\[
   \nbd{H}{D_1} \cap \nbd{h_igH}{D_1+\delta} \subseteq
     \nbd{h_iK}{D_2},
\]
where $D_2$ depends on $D_1$ but is independent of the choice of
$h_i \in H$.  In other words, $w$ lies within an $\Set{A}$--distance $D_2$
of both $h_1K$ and $h_2K$ in~$G$.  Since $H$ is an undistorted subgroup of~$G$,
it follows that the $\Set{B}$--distance between $h_1K$ and $h_2K$
in~$H$ is similarly bounded in terms of~$D$, as desired.
\end{proof}

\section{Relatively hyperbolic groups and the word metric}
\label{sec:RelHyp}

Various equivalent formulations of relatively hyperbolic groups have been
introduced and studied by Gromov \cite{Gromov87}, Farb \cite{Farb98},
Bowditch \cite{BowditchRelHyp},
\Drutu--Sapir \cite{DrutuSapirTreeGraded} and Osin \cite{Osin06}.

In this section we discuss geometric features of
the Cayley graph of a relatively hyperbolic group with respect to a finite
generating set.
In particular, we prove two results about triangles and quadrilaterals
with the property that each side lies near a peripheral subspace.

We begin by recalling the definition of a relatively hyperbolic group.
\begin{defn}[Relatively hyperbolic]
Let $G$ be a finitely generated
group and $\mathbb{P}$ a finite collection of subgroups of $G$.
Suppose $G$ acts on a $\delta$--hyperbolic graph $K$ with
finite edge stabilizers and finitely many orbits of edges.
Suppose $K$ is \emph{fine} in the sense that for each $n$
each edge of $K$ is contained in only finitely many circuits of length $n$.
Suppose also that $\mathbb{P}$ is a set of representatives of the
conjugacy classes of infinite vertex stabilizers.
Then the pair $(G,\mathbb{P})$ is \emph{relatively hyperbolic}.
The subgroups $P \in \mathbb{P}$ are the \emph{peripheral subgroups}
of $(G,\mathbb{P})$,
and their left cosets $gP$ are \emph{peripheral cosets}.
\end{defn}

Throughout this section all paths and distances are taken in the
Cayley graph of a relatively hyperbolic group with respect to a fixed
finite generating set.

We now collect several results due to \Drutu--Sapir
on the geometry of the word metric for a relatively hyperbolic group.
The first states that peripheral cosets are ``isolated.''

\begin{thm}[\Drutu--Sapir, Theorem~4.1, \cite{DrutuSapirTreeGraded}]
\label{thm:Isolated}
Suppose $(G,\mathbb{P})$ is relatively hyperbolic.
For each $\rho<\infty$ there is a constant $\kappa=\kappa(\rho)<\infty$
so that for any two peripheral cosets $gP \ne g'P'$ we have
\[
   \diam \bigl( \nbd{gP}{\rho} \cap \nbd{g'P'}{\rho} \bigr) < \kappa.
\]
\end{thm}

The next result quantifies the quasiconvexity of peripheral cosets
with respect to quasigeodesics.

\begin{thm}[\Drutu--Sapir, Lemma~4.15, \cite{DrutuSapirTreeGraded}]
\label{thm:QuasiconvexPeripheral}
Suppose $(G,\mathbb{P})$ is relatively hyperbolic.
Given positive constants $\epsilon$ and $\nu$, there is a constant
$\tau = \tau(\epsilon,\nu)$ so that the following holds.
Let $c$ be an $\epsilon$--quasigeodesic
with endpoints in the $\nu$--neighborhood of a peripheral coset $gP$.
Then $c \subseteq \nbd{gP}{\tau}$.
\end{thm}

The following proposition roughly states that
a pair of quasigeodesics beginning in the same peripheral coset
and ending in the same peripheral coset must leave the former coset and
enter the latter coset close together.

\begin{prop}[\Drutu--Sapir, Lemma~8.11, \cite{DrutuSapirTreeGraded}]
\label{prop:CloseEntryTwo}
Let $(G,\mathbb{P})$ be relatively hyperbolic.
There is a constant $\nu_0$ such that
for each $\nu > \nu_0$
there is a constant $D_1=D_1(\nu)$ such that the following holds.
Let $gP$ and $g'P'$ be two distinct peripheral cosets
and let $c_1$ and $c_2$ be two geodesics
such that $c_i$ has endpoints $x_i$ and~$y_i$.
If $c_i \cap \bar{\nbd{gP}{\nu}} = \{x_i\}$
and $c_i \cap \bar{\nbd{g'P'}{\nu}} = \{y_i\}$
then $d(x_1,x_2)$ and $d(y_1,y_2)$ are both less than $D_1$.
\end{prop}

\Drutu--Sapir also establish the following
geometric description of quasigeodesic triangles.

\begin{figure}[t]
\labellist
\small\hair 2pt
\pinlabel $c'_1$ at 112 64
\pinlabel $c'_2$ at 184 64
\pinlabel $c'_0$ at 142 7
\pinlabel $gP$ at 146 47
\endlabellist
\begin{center}
\includegraphics{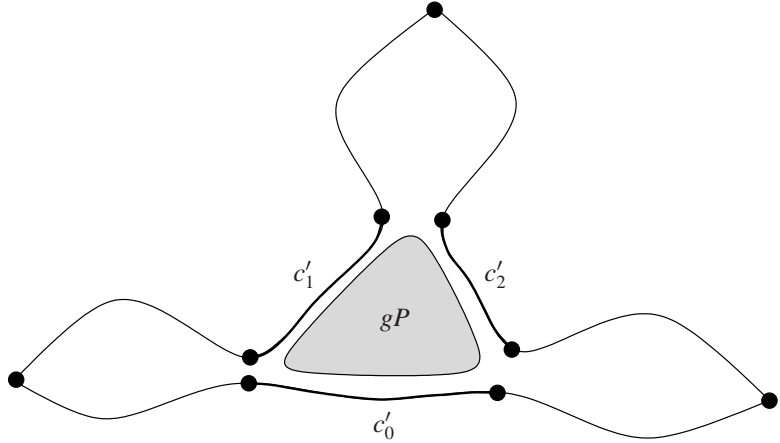}
\end{center}
\caption{A quasigeodesic triangle in a relatively hyperbolic group
with a peripheral coset at the center.}
\label{fig:Triangle}
\end{figure}

\begin{thm}[\Drutu--Sapir, Section 8.1.3, \cite{DrutuSapirTreeGraded}]
\label{thm:Triangle}
Let $(G,\mathbb{P})$ be relatively hyperbolic.
For each $\epsilon$ there are constants $\delta=\delta(\epsilon)$
and $\mu = \mu(\epsilon)$ such that the following holds.
Let $\Delta$ be an $\epsilon$--quasigeodesic triangle
with sides $c_0$, $c_1$ and $c_2$.
Then either
\begin{enumerate}
  \item there is a point $w$ such that the ball $\ball{w}{\delta}$
  intersects all three sides of~$\Delta$, or
  \item there is a peripheral coset $gP$ such that the neighborhood
  $\nbd{gP}{\delta}$ intersects all three sides of~$\Delta$.
\end{enumerate}
In the second case, illustrated in Figure~\ref{fig:Triangle},
let $c_i'$ be the smallest subpath of $c_i$
containing $c_i \cap \nbd{gP}{\delta}$.
Then
the terminal endpoint of $c_i'$
and the initial endpoint of $c'_{i+1}$
are mutually within a distance~$\mu$ \textup{(}indices modulo~$3$\textup{)}.
\end{thm}

Roughly speaking, the following lemma deals with a triangle each of whose sides
lies close to a peripheral coset.  The conclusion is that
either all three peripheral cosets are equal or one of the sides of the triangle
is short.

\begin{lem}
\label{lem:OneSideShort}
Let $(G,\mathbb{P})$ be relatively hyperbolic.
For each $\eta>0$, there is a constant $\zeta >0$ such that the following
holds.
Let $g_0 P_0$, $g_1 P_1$ and $g_2 P_2$ be peripheral cosets
such that $g_0 P_0 \notin \{g_1 P_1,g_2 P_2\}$.
Suppose for each $\{i,j,k\} = \{0,1,2\}$ there is a point
\[
   y_i \in \nbd{g_j P_j}{\eta} \cap \nbd{g_k P_k}{\eta}.
\]
Then
\[
   d(y_1,y_2) < \zeta.
\]
\end{lem}

\begin{proof}
For each choice of $\{i,j,k\} = \{0,1,2\}$,
pick a point $y_i \in \nbd{g_j P_j}{\eta} \cap \nbd{g_k P_k}{\eta}$.
Consider a geodesic triangle $\Delta:=c_0 \cup c_1 \cup c_2 \cup$
such that each $c_i$ has endpoints $y_j$ and $y_k$.
By Theorem~\ref{thm:QuasiconvexPeripheral}, there is a constant
$\tau$ depending on $\eta$ such that $c_i$ lies in the
$\tau$--neighborhood of $g_i P_i$.

Let $\delta = \delta(1)$ and $\mu=\mu(1)$ be the constants given by
Theorem~\ref{thm:Triangle}.
Then there is either a point $z$ or a peripheral coset $gP$
whose $\delta$--neighborhood intersects all three sides of $\Delta$.

We will see that we must be in the former case, for a possibly
larger constant $\delta'$.
Suppose $\nbd{gP}{\delta}$ intersects all three sides of $\Delta$.
For each $i$, let $c'_i$ denote the smallest subpath of $c_i$ containing
$c_i \cap \bar{\nbd{gP}{\delta}}$.

Choose $i \in \{0,1,2\}$ so that $gP \ne g_i P_i$.
Since $c'_i \subseteq \nbd{g_i P_i}{\delta}$, it follows from
Theorem~\ref{thm:Isolated} that $c'_i$ has length less than
$\kappa = \kappa(2\delta)$.
By Theorem~\ref{thm:Triangle}, the endpoints $z_j$ and $z_k$ of $c'_i$
are within a distance $\mu$ of the segments $c_j$ and $c_k$ respectively.
Thus we have
\[
   d(z_j,c_i) = 0, \quad d(z_j,c_j) < \mu \quad \text{and}
      \quad d(z_j,c_k) < \kappa + \mu.
\]
In other words, the ball of radius $\kappa+\mu$ centered at $z_j$
intersects all three sides of $\Delta$.

It follows that, in all cases, there exists a point $z$ such that the
$\delta'$--ball centered at $z$ intersects all three sides of $\Delta$,
where $\delta' = \delta + \kappa + \mu$.
Now for each choice of $\{j,k\} = \{1,2\}$ we have
\begin{gather*}
   z \in \nbd{g_0 P_0}{\delta' + \tau} \cap \nbd{g_j P_j}{\delta'+\tau}\\
\tag*{\text{and}}
   y_k \in \nbd{g_0 P_0}{\tau} \cap \nbd{g_j P_j}{\tau}.
\end{gather*}
Consequently, Theorem~\ref{thm:Isolated} gives an upper bound
$\kappa' := \kappa(\delta'+\tau)$
on the distance $d(y_k,z)$.
Varying $k \in \{1,2\}$, we see that
$d(y_1,y_2)$ is less than $2\kappa'$.
\end{proof}

The following result is analogous to the previous lemma, but for quadrilaterals instead of triangles.
If each side of a quadrilateral lies close to a peripheral
coset, then either all four cosets are equal or some pair of vertices
of the quadrilateral is close together.

\begin{lem}
\label{lem:FourSidesEqual}
Let $(G,\mathbb{P})$ be relatively hyperbolic.
For each $\eta>0$, there is a constant $\xi >0$ such that the following
holds.
Let $g_0 P_0$, $g_1 P_1$, $g_2 P_2$ and $g_3 P_3$ be peripheral cosets.
Suppose for each $i \in \{0,1,2,3\}$ there is a point
\[
   y_i \in \nbd{g_{i-1} P_{i-1}}{\eta} \cap \nbd{g_i P_i}{\eta}
\quad (\text{indices modulo }4).
\]
If the pairwise distances between the points $\{y_0,y_1,y_2,y_3\}$
are all at least $\xi$, then all four peripheral cosets are equal.
\end{lem}

\begin{proof}
Suppose the four peripheral cosets are not all equal.  We will show that
one of the pairwise distances between the $y_i$ is bounded above
by an appropriate constant $\xi$.

For each $i$, let $c_i$ be a geodesic in $X$
from $y_i$ to $y_{i+1}$.
Then by Theorem~\ref{thm:QuasiconvexPeripheral},
we have
\[
   c_i \subseteq \nbd{g_i P_i}{\tau}
\]
for some $\tau$ depending on $\eta$.

If $g_i P_i = g_{i+1} P_{i+1}$ for some $i$,
the result follows immediately from Lemma~\ref{lem:OneSideShort}.
Now suppose $g_i P_i=g_{i+2} P_{i+2}$ for some $i$, but
$g_{i+1} P_{i+1} \ne g_i P_i$.
Then $c_{i+1}$ lies in the $\tau$--neighborhood of both $g_i P_i$ and
$g_{i+1} P_{i+1}$.  But then Theorem~\ref{thm:Isolated}
gives an upper bound $\kappa_0 := \kappa(\tau)$
on $d(y_{i+1},y_{i+2})$, completing the proof.
Thus it suffices to assume that the four peripheral cosets are all
different.

By hypothesis, the points $y_0$ and $y_1$ lie in $\nbd{g_0 P_0}{\eta}$.
Let $u$ be the point where $c_3$ first enters the closed
$\eta$--neighborhood of $g_0 P_0$ when traversed from $y_3$ to $y_0$.
Similarly, let $v$ be the point where $c_1$ first enters
the closed $\eta$--neighborhood of $g_0 P_0$ when traversed from $y_2$
to $y_1$.
Since $g_0 P_0 \notin \{g_1 P_1,g_3 P_3\}$,
Theorem~\ref{thm:Isolated} gives an upper bound $\kappa_1:=\kappa(\eta)$
on the distances $d(y_0,u)$ and $d(y_1,v)$.
We may assume without loss of generality that $\eta$ is greater than the
constant $\nu_0$ given by Proposition~\ref{prop:CloseEntryTwo}.
Since $g_0 P_0 \ne g_2 P_2$,
Proposition~\ref{prop:CloseEntryTwo} gives an upper bound
$D_1=D_1(\eta)$ on the distance $d(u,v)$.
Thus we have
\[
   d(y_0,y_1) < D_1 + 2\kappa_1
\]
completing the proof of the lemma.
\end{proof}

\section{The geometry of relative Cayley graphs}
\label{sec:RelativeCayley}

The results in this section and the next
involve the interplay between two different
metrics on a relatively hyperbolic $G$ with a fixed finite generating
set $\Set{S}$.  When a metric is not specified, all distances are assumed
to refer to the word metric $d_{\Set{S}}$
or (equivalently) in the Cayley graph $\Cayley(G,d_{\Set{S}})$.
Let $\Set{P}$ be the union of all the peripheral
subgroups of $G$.
We will also use the word metric $d_{\Set{S} \cup \Set{P}}$ and the
Cayley graph $\Cayley(G,\Set{S} \cup \Set{P})$
for the (typically infinite) generating set $\Set{S} \cup \Set{P}$.
Whenever we use this ``relative metric'' or ``relative Cayley graph,''
we will explicitly call attention to its use.
In particular, if $A \subset G$ the notation
$\nbd{A}{\epsilon}$ always refers to
the $\epsilon$--neighborhood of $A$ 
using the metric $d_{\Set{S}}$,
and the notation $B(x,\epsilon)$ refers to an $\Set{S}$--metric ball.

Geodesics in the relative Cayley graph are sometimes referred to as
\emph{relative geodesics in $G$}, etc.
Observe that $\Cayley(G,\Set{S})$ is a subgraph of
$\Cayley(G,\Set{S} \cup \Set{P})$
containing all the vertices but omitting the edges labelled by elements
of $\Set{P}$.

The next result is a relatively hyperbolic analogue of the Morse Lemma,
dealing with a pair of relative geodesics whose respective
endpoints are close together in $d_{\Set{S}}$.

\begin{prop}[Osin, Proposition~3.15, \cite{Osin06}]
\label{prop:StrongGeodesicRelQC}
Let $(G,\mathbb{P})$ be relatively hyperbolic.
For each $k \ge 0$, there is a constant
$\rho=\rho(k)$ such that the following holds.
Let $c$ and $c'$ be relative geodesics in
$\Cayley(G,\Set{S} \cup \Set{P})$
such that their initial endpoints are within an $\Set{S}$--distance $k$
and their terminal endpoints are also within
an $\Set{S}$--distance $k$.
Then the set of vertices of $c$ and the set of vertices of $c'$
are within a Hausdorff $\Set{S}$--distance~$\rho$.
\end{prop}

The following result is an analogue of Theorem~\ref{thm:Triangle}
for relative geodesic triangles.

\begin{thm}[Osin, Theorem~3.26, \cite{Osin06}]
\label{thm:ThinTriangle}
Let $(G,\mathbb{P})$ be relatively hyperbolic.
Then there is a constant $\nu >0$ such that if $c_0 \cup c_1 \cup c_2$
is a relative geodesic triangle in $\Cayley(G,\Set{S} \cup \Set{P})$
then for each vertex $v$ of $c_0$ there is a vertex $u$
in the union $c_1 \cup c_2$
such that
\[
   d_{\Set{S}} (u,v) < \nu.
\]
\end{thm}

Subdividing an $n$--gon into $n-2$ triangles and applying the previous theorem,
we get the following corollary.  As an aside we note that
the linear function $n-2$ in the conclusion of the corollary
can be improved to a logarithmic function of $n$ by subdividing more carefully.

\begin{cor}
\label{cor:ThinNGon}
Let $(G,\mathbb{P})$ and $\nu$ be as in the previous theorem.
If $c_1 \cup \dots \cup c_n$ is a relative geodesic $n$--gon
in $\Cayley(G,\Set{S} \cup \Set{P})$,
then for each vertex $v$ on $c_1$ there is a vertex $u$
on the union $c_2 \cup \dots \cup c_n$ such that
\[
   d_{\Set{S}} (u,v) < (n-2)\nu.
\]
\end{cor}

The notion of saturation of a quasigeodesic was introduced
by \Drutu--Sapir \cite{DrutuSapirTreeGraded}.
The saturation of an arbitrary subset
is defined similarly, as follows.

\begin{defn}
Suppose $(G,\mathbb{P})$ is relatively hyperbolic,
and let $Y$ be a subset of~$G$.
For each $\nu>0$ the \emph{$\nu$--saturation} of~$Y$, denoted $\Sat_\nu(Y)$,
is the union of $Y$ and every peripheral coset intersecting
$\nbd{Y}{\nu}$.
\end{defn}

We conclude the section with a mild generalization of
Proposition~\ref{prop:StrongGeodesicRelQC}.

\begin{prop}
\label{prop:SaturationRelQC}
Let $(G,\mathbb{P})$ be relatively hyperbolic.
For each $\tau>0$ there is a constant $\lambda=\lambda(\tau)>0$ such that
the following holds.
Let $c$ and $c'$ be relative geodesic segments in
$\Cayley(G,\Set{S} \cup \Set{P})$,
such that the endpoints $x_0$ and $x_1$ of $c'$ lie within
an $\Set{S}$--distance $\tau$ of\/ $\Sat_\tau\bigl(\textup{Vert}(c)\bigr)$.
Let $v$ be a vertex of $c'$.
Then $v$ lies within an $\Set{S}$--distance $\lambda$
of either $x_0$, $x_1$, or a vertex of $c$.
\end{prop}

\begin{proof}
For each $i \in \{0,1\}$, we define a path $c_i$
in $\Cayley(G,\Set{S}\cup \Set{P})$ from $x_i$ to $c$ as follows.
If $x_i$ is within an $\Set{S}$--distance $\tau$ of a vertex of $c$, then
let $c_i$ be an $\Set{S}$--geodesic
of shortest length from $x_i$ to $c$.
Otherwise, $x_i$ is within an $\Set{S}$--distance $\tau$
of a left coset $g_iP_i$ that is within an $\Set{S}$--distance $\tau$ of
a vertex of $c$.
In the latter case, let $c_i$ be the concatenation of
an $\Set{S}$--geodesic $q_i$ of shortest length
from $x_i$ to $g_iP_i$
followed by a peripheral edge $e_i$ in $g_iP_i$ followed by
an $\Set{S}$--geodesic $r_i$ of shortest length from $g_iP_i$ to a vertex
of $c$.
Let $a_i\in C$ denote the terminal vertex of $c_i$,
and let $\bar{c}$ denote the portion of $c$ from $a_0$ to $a_1$.

Note that $c_i$ is a concatenation of at most $2\tau+1$ edges,
each of which is a relative geodesic of length one.
Thus $c' \cup c_0 \cup \bar{c} \cup c_1$ is
a relative geodesic polygon with at most $4\tau+4$ sides.
Consequently, by Corollary~\ref{cor:ThinNGon}
each vertex of $c'$ is within an $\Set{S}$--distance
$(4\tau+2)\nu$ of some vertex of $c_0 \cup \bar{c} \cup c_1$.

Choose a vertex $v$ of $c'$ such that $v$ is at an $\Set{S}$--distance
more than $(4\tau+2)\nu + \tau$ from $x_0$.
If $c_0$ is an $\Set{S}$--geodesic of length at most $\tau$,
then, by the triangle inequality, no vertex of $c_0$
can be within an $\Set{S}$--distance $(4\tau+2)\nu$ of $v$.
Thus $v$ is within an $\Set{S}$--distance $(4\tau+2)\nu$ of
some vertex of $\bar{c} \cup c_1$.
On the other hand, if $c_0$ is a concatenation $q_0 \cup e_0 \cup r_0$
as described above,
then $q_0$ has $\Set{S}$--length at most $\tau$.
By a similar argument, it follows that $v$ is within an $\Set{S}$--distance
$(4\tau+2)\nu$ of some vertex of $r_0 \cup \bar{c} \cup c_1$.
Since $r_0$ has $\Set{S}$--length at most $\tau$, each of its vertices lies
within an $\Set{S}$--distance $\tau$ of the vertex $a_0$.
Thus in either case, $v$ is within an $\Set{S}$--distance
$(4\tau+2)\nu + \tau$ of some vertex of $\bar{c} \cup c_1$.

Interchanging the roles of $c_0$ and $c_1$, we see that
if $v$ is also at least an $\Set{S}$--distance $(4\tau+2)\nu + \tau$ from $x_1$
then it is within an $\Set{S}$--distance $(4\tau+2)\nu +\tau$ from some
vertex of $\bar{c}$.
Setting $\lambda = (4\tau+2)\nu + \tau$ completes the proof.
\end{proof}

\begin{defn}
Let $(G,\mathbb{P})$ be relatively hyperbolic with a finite generating set
$\Set{S}$.
Let $c$ be a geodesic of $\Cayley(G,\Set{S})$
and let $\epsilon$ and $R$ be positive constants.
A point $x \in c$ is \emph{$(\epsilon,R)$--deep} in a peripheral coset $gP$
(with respect to $c$) if $x$ is not within an $\Set{S}$--distance $R$
of an endpoint of
$c$ and if $B(x,R) \cap c$ lies in $\nbd{gP}{\epsilon}$.
If $x$ is not $(\epsilon,R)$--deep in any $gP$
then $x$ is an \emph{$(\epsilon,R)$--transition point} of $c$.
\end{defn}

\begin{prop}[Hruska, Proposition~8.13, \cite{HruskaRelQC}]
\label{prop:TransitionPoints}
Let $(G,\mathbb{P})$
be a relatively hyperbolic group with a finite generating set
$\Set{S}$.
There exist constants $\epsilon$, $R$ and $L$ such that the following holds.
Let $c$ be any geodesic of $\Cayley(G,\Set{S})$ with endpoints in $G$,
and let $c'$ be a relative geodesic in $\Cayley(G,\Set{S} \cup \Set{P})$
with the same endpoints as $c$.
Then the set of vertices of $c'$ is at a Hausdorff
$\Set{S}$--distance at most $L$
from the set of $(\epsilon,R)$--transition points of $c$.

Furthermore, let $\bar{c}$ be a connected component of the set of
all $(\epsilon,R)$--deep points of $c$.
Then there is a peripheral coset $gP$ such that each $x \in \bar{c}$
is $(\epsilon,R)$--deep in $gP$ and is not $(\epsilon,R)$--deep in
any other peripheral coset.
\end{prop}

\section{Relatively quasiconvex subspaces}
\label{sec:RelQC}

In the theory of word hyperbolic groups, the most natural subgroups
(and subspaces) are the quasiconvex subgroups (and subspaces).
In this section, we examine ``relatively quasiconvex'' subspaces,
which play an analogous role in the theory of relatively hyperbolic
groups.
The definition of relative quasiconvexity given below
was introduced by Osin \cite{Osin06}.

\begin{defn}
Let $(G,\mathbb{P})$ be relatively hyperbolic such that $G$ is generated by
a finite set $\Set{S}$.
For $\sigma\geq 0$, a subset $A \le G$ is
\emph{$\sigma$--relatively quasiconvex} in $G$
if the following condition holds:
Let $c$ be a relative geodesic in
$\Cayley(G,\Set{S} \cup \Set{P})$
connecting two points of $Y$.
Then each vertex $v$ of $c$ lies in the $\sigma$--neighborhood of $A$.

A subset $A \subseteq G$ is \emph{relatively quasiconvex}
if it is $\sigma$--relatively quasiconvex for some $\sigma \geq 0$.
\end{defn}

\begin{prop}[Osin, Proposition~4.10, \cite{Osin06}]
Relative quasiconvexity of a subset $A \subseteq G$
is independent of the choice of finite generating set $\Set{S}$ for $G$.
\end{prop}

The following result analogous to Theorem~\ref{thm:Triangle}
follows immediately from the definition of relative quasiconvexity.

\begin{prop}\label{prop:PeripheralRelQC}
Each peripheral coset $gP$ in a relatively hyperbolic group is
$0$--relatively quasiconvex.\qed
\end{prop}

The following result was proved independently by
Hruska \cite[Theorem~9.1]{HruskaRelQC}
and Martinez-Pedroza \cite[Proposition~1.5]{MartinezPedroza}.

\begin{prop}\label{prop:FinConjClasses}
Let $(G,\mathbb{P})$ be relatively hyperbolic and let $H \le G$
be $\sigma$--relatively quasiconvex.
The infinite maximal parabolic subgroups $Q=H\cap gPg^{-1}$
of $H$ lie in finitely many $H$--conjugacy classes,
and the corresponding peripheral cosets $gP$ in $G$
lie in finitely many $H$--orbits.
\end{prop}

\begin{prop}
\label{prop:Center}
Let $(G,\mathbb{P})$ be relatively hyperbolic.
Choose positive constants $\sigma$ and~$\tau$.
Then there is a constant $R = R(\sigma,\tau)$ such that the following holds.
Choose $n \ge 3$, and let $\Set{A} = \{A_1,\dots,A_n\}$
be a collection of $\sigma$--relatively quasiconvex subsets of~$G$.
Suppose for each $i = 1,\dots,n$ there is a point
\[
   x_i \in \bigcap
     \bigset{\nbd{A_j}{\tau}}{1\leq j \leq n,\ j\neq i}.
\]
Then we have the following consequences:
\begin{enumerate}
\item
\label{item:Center2}
There exists a vertex
\[
   x_0 \in \bignbd{\Sat_{R}(A_1)}{R} \cap \nbd{A_2}{R} \cap \dots
     \cap \nbd{A_n}{R}.
\]
\item
\label{item:Center2a}
If $A_1$ is a left coset of a peripheral subgroup, then there exists a
vertex 
\[
   x_0 \in \nbd{A_1}{R} \cap \dots \cap \nbd{A_n}{R}.
\]
\item
\label{item:Center3}
If $A_0$ is another $\sigma$--relatively quasiconvex subset of~$G$
such that $x_i$ is contained in
$\bignbd{\Sat_\tau (A_0)}{\tau}$ for all $i=1,\dots,n$,
then there exists a vertex $x_0$ satisfying
(\ref{item:Center2a})
such that $x_0 \in \bignbd{\Sat_R (A_0)}{R}$.
\end{enumerate}
\end{prop}

\begin{proof}
Choose relative geodesics
\[
   c_1 = [x_2,x_3], \quad c_2 = [x_3,x_1] \quad \text{and}
   \quad c_3 = [x_1,x_2]
\]
in $\Cayley(G,\Set{S} \cup \Set{P})$,
and consider the triangle $\Delta := c_1 \cup c_2 \cup c_3$.
By Proposition~\ref{prop:StrongGeodesicRelQC} and the $\sigma$--relative
quasiconvexity of $A_\ell$, there is a constant $\rho=\rho(\tau)$
such that if $i \in \{1,2,3\}$ and $\ell \in \{i\} \cup \{4,\dots,n\}$
then each vertex of $c_i$ lies within an $\Set{S}$--distance $\sigma+\rho$
of $A_\ell$.

By Theorem~\ref{thm:ThinTriangle}, there is a constant $\nu$
such that for each choice of $\{i,j,k\} = \{1,2,3\}$
each vertex of $c_i$ lies within an $\Set{S}$--distance $\nu$
of a vertex of $c_j \cup c_k$.
It follows that $c_i$ contains a pair of adjacent vertices $v_{ij}$
and $v_{ik}$
connected by an edge $e_i$ such that $v_{ij}$ is within an $\Set{S}$--distance
$\nu$ of a vertex of $c_j$ and $v_{ik}$ is within an $\Set{S}$--distance
$\nu$ of a vertex of $c_k$.

If for some choice of $\{i,j,k\} = \{1,2,3\}$,
the edge $e_i$ is labelled by a generator in $\Set{S}$, then
for all $\ell \in \{i\} \cup \{4, \dots, n\}$ we have:
\begin{equation}
\label{eqn:TripleCenter}
   v_{ij} \in \nbd{A_k}{1 + \nu + \sigma + \rho}
             \cap \nbd{A_j}{\nu + \sigma + \rho}
         \cap \nbd{A_\ell}{\sigma + \rho}
\end{equation}

If not, then $e_3$ is labelled by a generator in $\Set{P}$
and $v_{32}$ and $v_{31}$ are vertices of a peripheral coset $gP$.
Since $d_{\Set{S}}(v_{31},A_1) < \nu + \sigma + \rho$, the coset
$gP$ is contained in $\Sat_{\nu+\sigma + \rho}(A_1)$.
Whenever $\ell=3,\dots,n$, we now have:
\[
   v_{32} \in \Sat_{\nu + \sigma + \rho}(A_1)
            \cap \nbd{A_2}{\nu + \sigma + \rho}
        \cap \nbd{A_\ell}{\sigma + \rho}
\]
In either case, choosing $x_0 := v_{32}$ and $R \ge 1 + \nu + \sigma + \rho$
completes the proof of (\ref{item:Center2}).

Now let us consider assertion (\ref{item:Center2a}).
Suppose the subspace $A_1$ is a peripheral coset $gP$.
If for some $i \in \{1,2,3\}$ the edge $e_i$ corresponds
to a generator in $\Set{S}$ then we are done,
since (\ref{eqn:TripleCenter}) shows that
$v_{ij}$ is within an $\Set{S}$--distance
$1 + \nu + \sigma + \rho$ of $A_\ell$ for all $\ell = 1,\dots,n$.

It therefore suffices to assume that each $e_i$ corresponds to a
generator from $\Set{P}$ such that $v_{ij},v_{ik} \in g_i P_i$.
for some peripheral coset $g_iP_i$.
If $gP \ne g_1P_1$, then by Theorem~\ref{thm:Isolated},
the $\Set{S}$--distance between $v_{12}$ and $v_{13}$ is at most
$\kappa = \kappa(\sigma + \rho)$.
In this case, a computation similar to (\ref{eqn:TripleCenter})
shows that $v_{12}$ is within an $\Set{S}$--distance
$\kappa + \nu + \sigma + \rho$ of $A_\ell$ for all $\ell = 1,\dots,n$.
Similarly
if the cosets $g_1P_1$, $g_2P_2$ and $g_3P_3$ are not all the same,
then by Lemma~\ref{lem:OneSideShort}
the $\Set{S}$--distance between some pair $v_{ij}$ and $v_{ik}$
is bounded above by $\zeta = \zeta(\nu)$.
Therefore $v_{ij}$ is within an $\Set{S}$--distance
$\zeta + \nu + \sigma + \rho$ of $A_\ell$ for all $\ell = 1,\dots,n$.

The only remaining possibility is that the four cosets
$gP$, $g_1P_1$, $g_2P_2$ and $g_3P_3$ are all equal.
In particular, we have $A_1 = gP = g_3P_3$.
It follows that
\[
   d_{\Set{S}}(v_{32},A_1) = 0, \quad
   d_{\Set{S}}(v_{32},A_2) < \nu + \sigma + \rho \quad \text{and} \quad
   d_{\Set{S}}(v_{32},A_\ell) < \sigma + \rho
\]
for all $\ell = 3,\dots,n$.
In all possible cases, to complete the proof of (\ref{item:Center2a})
it suffices to choose
\[
   R > 1 + \kappa + \zeta + \nu + \sigma + \rho.
\]

Finally, we turn our attention to assertion (\ref{item:Center3}).
We have shown that assertions
(\ref{item:Center2}) and~(\ref{item:Center2a})
hold for sufficiently large $R$ and for some vertex $x_0$ lying
on one of the sides $c_i$ of the triangle $\Delta$.
Suppose the endpoints $x_j$ and $x_k$ of $c_i$
lie within an $\Set{S}$--distance $\tau$
of $\Sat_\tau(A_0)$.
Choose points $u$ and $v$ in $A_0$ such that
$x_j \in \bignbd{\Sat_\tau(u)}{\tau}$ and
$x_k \in \bignbd{\Sat_\tau(v)}{\tau}$,
and let $\bar{c}$ be a geodesic in $\Cayley(G,\Set{S} \cup \Set{P})$
from $u$ to $v$.

It follows from Proposition~\ref{prop:SaturationRelQC}
that the vertex $x_0$ of $c_i$ lies within a distance
$\lambda=\lambda(\tau)$ of either $x_j$, $x_k$, or a vertex of $\bar{c}$.
Since $A_0$ is $\sigma$--quasiconvex, it follows that $x_0$
lies within a distance $\lambda + \sigma$ of $\{x_j,x_k\} \cup A_0$.
In particular, $x_0$ is within an $\Set{S}$--distance $\lambda + \sigma + \tau$
of $\Sat_\tau(A_0)$, establishing (\ref{item:Center3})
for all $R > \lambda + \sigma + \tau$.
\end{proof}

\begin{cor}
\label{cor:AllNearPeripheral}
Let $(G,\mathbb{P})$ be relatively hyperbolic.
Choose positive constants $\sigma$ and $D$ and an integer $n\ge 1$.
There is a constant $R_1 = R_1(\sigma,D,n)$ such that the following holds.
Let $\Set{A}=\{A_1,\dots,A_n\}$ be a set of $\sigma$--relatively quasiconvex
subsets of $G$.  Suppose the following intersection is nonempty:
\[
   \bigcap_{i=1}^n \nbd{A_i}{D}
\]
If there is a peripheral coset $gP$ such that for all $i$
the intersection $\nbd{A_i}{D} \cap \nbd{gP}{D}$ is nonempty,
then the following intersection is nonempty:
\[
   \nbd{gP}{R_1} \cap \nbd{A_1}{R_1} \cap \cdots \cap \nbd{A_n}{R_1}
\]
\end{cor}

\begin{proof}
Recall that $gP$ is $0$--relatively hyperbolic.
When $n=1$, the result is trivial.
We proceed by induction on $n$ for $n\ge 2$.
Let $\tau:= R_1(\sigma,D,n-1)$ denote the constant given by the
inductive hypothesis.
By assumption, $\bigcap_{i=1}^n \nbd{A_i}{D}$
is nonempty.  By induction, for each $i=1,\dots,n$ the following
intersection is nonempty:
\[
   \nbd{gP}{\tau} \cap
      \bigcap \bigset{\nbd{A_j}{\tau}}{1\le j\le n,\ j\neq i}
\]
The result now follows from
Proposition~\ref{prop:Center}(\ref{item:Center2a}),
applied to the collection $\{gP,A_1,\dots,A_n\}$
of $\sigma$--quasiconvex subspaces.
\end{proof}

\begin{prop}\label{prop:CommonIntersection}
Let $(G,\mathbb{P})$ be relatively hyperbolic.
Choose positive constants $\sigma$ and $D$ and an integer $n \ge 1$.
Let $\Set{A}$
be an arbitrary set of $\sigma$--relatively quasiconvex subsets of~$G$
such that
for each $A,A' \in \Set{A}$ we have $d_{\Set{S}}(A,A') < D$.
Then there is a constant $M_0=M_0(\sigma,D,n)$ so that whenever $M\ge M_0$
at least one of the following holds.
\begin{enumerate}
  \item For every subset $\{A_1,\dots,A_n\}\subseteq\Set{A}$, the following intersection is nonempty:
  \[
     \bigcap_{\ell=1}^{n} \nbd{A_\ell}{M}
  \]
  \item
  \label{item:CI:Peripheral}
  There is a peripheral coset $gP$ such that $d_{\Set{S}}(gP,A) < M$
  for all $A \in \Set{A}$.
  Furthermore, for all $A,A' \in \Set{A}$, the following triple intersection
  is nonempty:
\[
   \nbd{gP}{M} \cap \nbd{A}{M} \cap \nbd{A'}{M}
\]
\end{enumerate}
\end{prop}

\begin{proof}
Recall that every peripheral coset $gP$ is $0$--quasiconvex,
as observed in Proposition~\ref{prop:PeripheralRelQC}.
In order to establish (\ref{item:CI:Peripheral}),
it is sufficient to find a constant $M$ and a peripheral coset $gP$
such that $d_{\Set{S}}(gP,A)<M$
for all $A \in \Set{A}$.
For then applying Proposition~\ref{prop:Center}(\ref{item:Center2a})
to each collection $\{gP,A,A'\}$ with $A,A' \in \Set{A}$
gives a larger constant $M'$ so that the following triple intersection
is nonempty:
\[
   \nbd{gP}{M'} \cap \nbd{A}{M'} \cap \nbd{A'}{M'}
\]

The proposition is trivial when $n\leq 2$.
We proceed by induction on~$n$ for $n\ge 3$.
Let $\tau$ denote the constant $M_0(\sigma,D,n-1)$ given by the inductive
hypothesis.
We will show that the proposition holds whenever $M$
is sufficiently large (in terms of $\sigma$ and $\tau$).

Suppose $\{A_0,\dots,A_{n-1}\}$ is a subset of~$\Set{A}$ such that
\begin{equation}
\label{eqn:Empty1}
   \bigcap_{\ell=0}^{n-1} \nbd{A_\ell}{M} = \emptyset.
\end{equation}
If no such set exists, there is nothing to prove.
Furthermore, if there is a peripheral coset $gP$
such that $d_{\Set{S}}(P,A) < \nu$ for all $A \in \Set{A}$, then we are done.
Thus we may assume
that for every collection of $n-1$ distinct elements
of $\Set{A}$, their $\nu$--neighborhoods
have a nonempty intersection.

Applying Proposition~\ref{prop:Center}(\ref{item:Center2})
to the collection $\{A_0,\dots,A_{n-1}\}$ gives $R_0=R_0(\sigma,\tau)$
and a point
\[
   x_n \in \bignbd{ \Sat_{R_0} (A_0) }{R_0}
            \cap \nbd{A_1}{R_0} \cap \dots \cap \nbd{A_{n-1}}{R_0}.
\]
If we choose $M$ greater than $R_0$, then (\ref{eqn:Empty1})
implies that $d_{\Set{S}}(x_n,A_0)$ is at least~$R_0$.
Thus $x_n \in \nbd{g_0P_0}{R_0}$ for some peripheral coset $g_0P_0$
such that $g_0P_0 \subseteq \Sat_{R_0}(A_0)$.
Note that $d_{\Set{S}}(g_0P_0,A_\ell) < 2R_0$ for all $\ell = 0,\dots,n-1$.
If $\Set{A} = \{A_0,\dots,A_{n-1}\}$, we are done.
Otherwise, choose an arbitrary
$A_n \in \Set{A} \setminus \{A_0,\dots,A_{n-1}\}$.
To complete the proof,
we will bound the $\Set{S}$--distance from $A_n$ to $g_0P_0$
in terms of $\sigma$ and $\tau$, provided that $M$ is sufficiently large.

For each $i = 1,\dots,n-1$,
applying Proposition~\ref{prop:Center}(\ref{item:Center2})
to the collection
\[
   \{A_0,\dots,A_{i-1},A_{i+1},\dots,A_n\}
\]
gives a point
\[
   x_i \in \bignbd{\Sat_{R_0}(A_0)}{R_0}
             \cap \nbd{A_1}{R_0} \cap \dots
         \cap \nbd{A_{i-1}}{R_0} 
         \cap \nbd{A_{i+1}}{R_0} \cap \dots \cap \nbd{A_n}{R_0}.
\]
Notice that the points $x_1,\dots,x_n$ and subspaces $A_0,A_1,\dots,A_n$
satisfy the hypothesis of Proposition~\ref{prop:Center}(\ref{item:Center3}).
Therefore there exist a constant $R_1:=R(\sigma,R_0)$
and a point
\[
   x_0 \in \bignbd{\Sat_{R_1}(A_0)}{R_1}
             \cap \nbd{A_1}{R_1} \cap \dots \cap \nbd{A_n}{R_1}.
\]

Choose a relative geodesic $c$ from $x_0$ to $x_n$.
Since the endpoints of $c$ lie within a distance $R_0+R_1$ of $A_\ell$
for each $\ell=1,\dots,n-1$, we know from
Proposition~\ref{prop:StrongGeodesicRelQC} and the $\sigma$--relative
quasiconvexity of $A_\ell$ that each vertex of $c$ lies within an
$\Set{S}$--distance
$\sigma + \rho$ of $A_\ell$, where $\rho = \rho(R_0+R_1)$.
However, by (\ref{eqn:Empty1}) there is no vertex within a distance
$M$ of every $A_\ell$ for $\ell=0,\dots,n$.

Roughly speaking, it follows that once $M$ is sufficiently large,
each vertex of $c$ will be far from $A_0$.
We have already used this principle to conclude that
$x_n$ is within an $\Set{S}$--distance $R_0$ of $g_0P_0$,
which is itself within
an $\Set{S}$--distance $R_0$ of $A_0$, once $M$ is sufficiently large.
Let $y_n$ be a vertex of $A_0$ within an $\Set{S}$--distance $R_0$
of $g_0P_0$.
A similar argument shows that $x_0$ is within a distance $R_1$
of a peripheral coset $g_1P_1$,
which is itself within a distance $R_1$ of $A_0$,
for sufficiently large $M$.
Let $y_0$ be a vertex of $A_0$ within an $\Set{S}$--distance $R_1$ of $g_1P_1$.
Choose a relative geodesic $\bar{c}$ of $\Cayley(G,\Set{S} \cup \Set{P})$
from $y_0$ to $y_n$.

Since $x_0$ and $x_n$ are both contained in
$\bignbd{ \Sat_{R_0 + R_1}(\bar{c})}{R_0 + R_1}$,
it follows from Proposition~\ref{prop:SaturationRelQC}
that each vertex of $c$ lies within an $\Set{S}$--distance $\lambda$
of either $x_0$, $x_n$, or a vertex of $\bar{c}$.
By the $\sigma$--relative quasiconvexity of $A_0$,
each vertex of $\bar{c}$ lies within an $\Set{S}$--distance $\sigma$ of $A_0$.
Thus each vertex of $c$ either lies within an $\Set{S}$--distance $\lambda$
of $\{x_0,x_n\}$
or lies within an $\Set{S}$--distance $\lambda + \sigma$ of $A_0$ for some
$\lambda = \lambda(R_0 + R_1)$.
Suppose $v$ is a vertex of $c$ whose $\Set{S}$--distance from $A_0$
is less than $\lambda + \sigma$.
Once $M > \lambda + \sigma + \rho$, the existence of such a vertex $v$
contradicts (\ref{eqn:Empty1}).
Therefore we may assume that each vertex $v$ of $c$
lies within an $\Set{S}$--distance $\lambda$ of $\{x_0,x_n\}$.
A similar argument shows that each vertex $v$ of $\bar{c}$
lies within an $\Set{S}$--distance $\lambda$ of $\{y_0,y_n\}$,
since otherwise $v$ would be close to both $A_0$ and a vertex of $c$,
contradicting (\ref{eqn:Empty1}) once $M$ is sufficiently large.

Let $e$ be an edge of $c$ with endpoints $v_0$ and $v_n$
such that $v_0$ is within an $\Set{S}$--distance $\lambda$ of $x_0$
and $v_n$ is within an $\Set{S}$--distance $\lambda$ of $x_n$.
If $e$ corresponds to a generator in $\Set{S}$, then
$d_{\Set{S}}(x_0,x_n) < 2\lambda + 1$.
Since $d_{\Set{S}}(x_n,g_0P_0) < R_0$ and $d_{\Set{S}}(x_0,A_n) < R_1$,
it follows that
$d_{\Set{S}}(A_n,g_0P_0) < 2\lambda + 1 + R_0 + R_1$ and we are done.
Thus we may assume that $e$ corresponds to a generator in $\Set{P}$
and that $v_0,v_n \in g_2P_2$ for some peripheral coset $g_2P_2$.

Now let $\bar{e}$ be an edge of $\bar{c}$ with endpoints $w_0$ and $w_n$
within a distance $\lambda$ of $y_0$ and $y_n$ respectively.
We will see that $\bar{e}$ must also correspond to a generator in
$\Set{P}$.
Suppose instead that $\bar{e}$ corresponds to a generator in $\Set{S}$.
Then $d_{\Set{S}}(y_0,y_n)$ is less than $2\lambda + 1$.
In this case, consider the three peripheral cosets $g_0P_0$, $g_1P_1$,
and $g_2P_2$.
Recall that $x_0$ lies within a distance $R_1 + \lambda$ of both
$g_1P_1$ and $g_2P_2$.
If $g_0P_0 \in \{g_1P_1,g_2P_2\}$, then $d(x_0,g_0P_0) < R_1 + \lambda$,
completing the proof as before since $d(x_0,A_n) < R_0$.
On the other hand, if $g_0P_0 \notin \{g_1P_1,g_2P_2\}$,
we can apply Lemma~\ref{lem:OneSideShort} to $g_0P_0$, $g_1P_1$ and $g_2P_2$
to produce an upper bound $\zeta$ on the distance from $x_n$ to $A_0$,
contradicting (\ref{eqn:Empty1}) once $M$ is at least
$\zeta + \sigma + \rho$ as before.
Therefore we may assume that $\bar{e}$ corresponds to a generator
in $\Set{P}$ and that $w_0,w_n \in g_3P_3$
for some peripheral coset $g_3P_3$.

We now have a cycle of four peripheral cosets $g_0P_0$, $g_2P_2$, $g_1P_1$,
and $g_3P_3$ such that the following hold for $\eta := R_0 + R_1 + \lambda$:
\begin{align*}
   x_n &\in \nbd{g_0P_0}{\eta} \cap \nbd{g_2P_2}{\eta} \\
   x_0 &\in \nbd{g_2P_2}{\eta} \cap \nbd{g_1P_1}{\eta} \\
   y_0 &\in \nbd{g_1P_1}{\eta} \cap \nbd{g_3P_3}{\eta} \quad \text{and} \\
   y_n &\in \nbd{g_3P_3}{\eta} \cap \nbd{g_0P_0}{\eta}
\end{align*}
Now Lemma~\ref{lem:FourSidesEqual} implies that either all four cosets are
equal, or some pair from the set $\{x_n,x_0,y_0,y_n\}$ has distance
less than $\xi$, for some constant $\xi = \xi(\eta)$.

We can now complete the proof using arguments similar to the ones used above.
As we have seen, it suffices to assume that
$g_0P_0 \notin \{g_1P_1,g_2,P_2\}$.
Therefore some pair from $\{x_n,x_0,y_0,y_n\}$ is within a distance $\xi$.
If $d(y_0,y_n) < \xi$, then the result follows from an application of
Lemma~\ref{lem:OneSideShort}.
If $d(x_0,x_n) < \xi$, then $d(A_n,g_0P_0) <  R_0 + R_1 + \xi$.
If $d\bigl( \{x_0,x_n\}, \{y_0,y_n\} \bigr) < \xi$
then, since $y_0,y_n \in A_0$, it follows that
either $x_0$ or $x_n$ is within a distance $\xi$ of $A_0$
and within a distance $R_0 + R_1$ of $A_1,\dots,A_{n-1}$,
contradicting (\ref{eqn:Empty1}) once $M > \xi + R_0 + R_1$.
\end{proof}

\begin{prop}[Close to two peripherals]
\label{prop:NearTwoPeripherals}
Let $(G,\mathbb{P})$ be relatively hyperbolic.
For each $M>0$ and $\sigma>0$
there exists $M_1=M_1(M,\sigma)$ such that the following holds.
Suppose $\Set{A}$ is a collection of $\sigma$--relatively quasiconvex
subsets and $g_0P_0\ne g_1P_1$ are distinct peripheral cosets.
If $d_{\Set{S}}(g_iP_i,A)<M$ for all $i\in \{0,1\}$
and for all $A \in \Set{A}$,
then the following intersection is nonempty:
\[
   \bigcap_{A\in\Set{A}} \nbd{A}{M_1}
\]
\end{prop}

\begin{proof}
Let $\epsilon$, $R$ and $L$ be the constants given by
Proposition~\ref{prop:TransitionPoints},
let $\nu_0$ be the constant given by
Proposition~\ref{prop:CloseEntryTwo},
let $\nu := \nu_0 + M + \epsilon$,
and let $D_1=D_1(\nu)$ be given by Proposition~\ref{prop:CloseEntryTwo}.
Let $\tau=\tau(\nu)$ be given by Theorem~\ref{thm:QuasiconvexPeripheral},
and let $\kappa=\kappa(\epsilon+\tau)$ be given by
Theorem~\ref{thm:Isolated}.

For each $A \in \Set{A}$ there exists a geodesic $c$
in $\Cayley(G,\Set{S})$
whose endpoints $x_0,x_1$ satisfy $d_{\Set{S}}(x_i,g_iP_i) < M$.
Choose $y \in c$ to be the last point of $c$ that lies in
$\overline{\nbd{g_0P_0}{\nu}}$
(when $c$ is traversed from $x_0$ to $x_1$).

We will see that $y$ is within an $\Set{S}$--distance $\kappa$ of
an $(\epsilon,R)$--transition point
of $c$.  Suppose $y$ is not itself an $(\epsilon,R)$--transition point.
Then by Proposition~\ref{prop:TransitionPoints}
there is a subgeodesic $\bar{c}$ of $c$ containing $y$
each of whose interior points
is $(\epsilon,R)$--deep in a fixed peripheral coset $gP$
and whose endpoints are $(\epsilon,R)$--transition points.
It is clear that $gP \ne g_0P_0$, since $\epsilon \le \nu$.
Let $[x_0,y]$ denote the portion of $c$ from $x_0$ to $y$.
Since $M \le \nu$,
the points $x_0$ and $y$ each lie in $\nbd{g_0P_0}{\nu}$.
Therefore, $[x_0,y]$ lies in $\nbd{g_0P_0}{\tau}$.
On the other hand, $\bar{c}$ lies in $\nbd{gP}{\epsilon}$.
Thus $\bar{c} \cap [x_0,y]$ lies in $\nbd{gP}{\epsilon} \cap
\nbd{g_0P_0}{\tau}$ and hence has length at most $\kappa$.
In other words, $y$ is within an $\Set{S}$--distance $\kappa$
of an $(\epsilon,R)$--transition point $z$ of $c$, as desired.

By Proposition~\ref{prop:TransitionPoints}
this point $z$ is within an $\Set{S}$--distance $L$
of a vertex $w$ of $\hat{c}$, where $\hat{c}$ is any relative geodesic
in $\Cayley(G,\Set{S} \cup \Set{P})$ with the same endpoints as $c$.
Since $A$ is $\sigma$--relatively quasiconvex, it follows that
$d_{\Set{S}}(w,A) < \sigma$.  Therefore,
$d_{\Set{S}}(y,A) < \kappa + L + \sigma$.

If we now vary the choice of $A$, the above argument
produces for each $A \in \Set{A}$ a point $y = y(A)$
within an $\Set{S}$--distance $\kappa + L + \sigma$ of $A$.
Since $y$ is the last point of $c$ that lies in
$\overline{\nbd{g_0P_0}{\nu}}$,
it follows from Proposition~\ref{prop:CloseEntryTwo} that
for each $A,A' \in \Set{A}$ we have
$d_{\Set{S}}\bigl( y(A),y(A') \bigr) < D_1$.
In particular, if we fix $A \in \Set{A}$,
then every $A' \in \Set{A}$ comes within an $\Set{S}$--distance
$D_1 + \kappa + L + \sigma$ of the point $y(A)$.
Setting $M_1 \ge D_1 + \kappa + L + \sigma$ completes the proof.
\end{proof}

The following corollary essentially states that, if we increase the constants given by Proposition~\ref{prop:CommonIntersection}, then the peripheral
coset appearing in conclusion (\ref{item:CI:Peripheral}) is unique.

\begin{cor}[Uniqueness of peripheral coset]
\label{cor:UniquePeripheral}
Let $(G,\mathbb{P})$ be relatively hyperbolic.
Choose $\sigma$, $D$, $n$ and $\Set{A}$ satisfying the hypothesis of
Proposition~\ref{prop:CommonIntersection},
and let $M_0 = M_0(\sigma,D,n)$ be the constant given by
Proposition~\ref{prop:CommonIntersection}.
Then for each $M \ge M_0$ there exists $M_2=M_2(\sigma,D,n,M)$
such that at least one of the following holds:
\begin{enumerate}
\item
\label{item:Unique:Center}
For every $\{A_1,\dots,A_n\} \subseteq \mathcal{A}$, the 
following intersection is nonempty:
\[
   \bigcap_{\ell=1}^n \nbd{A_\ell}{M_2}
\]
\item
\label{item:Unique:Peripheral}
There is a unique peripheral coset $gP$ such that
$d_{\Set{S}}(gP,A) < M$ for all $A \in \Set{A}$.
Furthermore, for all $A,A' \in \Set{A}$ the following intersection
is nonempty:
\[
   \nbd{gP}{M} \cap \nbd{A}{M} \cap \nbd{A'}{M}
\]
\end{enumerate}
\end{cor}

\begin{proof}
Choose $M\ge M_0$ arbitrary.
If there does not exist a peripheral coset $gP$ such that
for all $A \in \Set{A}$ we have $d_{\Set{S}}(gP,A)< M$, then
Proposition~\ref{prop:CommonIntersection} implies that (\ref{item:Unique:Center}) holds for any $M_2 \ge M$.

On the other hand, suppose there are two
peripheral cosets $g_0P_0 \ne g_1P_1$
such that $d_{\Set{S}}(g_i P_i, A) < M$ for all $i\in\{0,1\}$ and all
$A \in \Set{A}$.
In this case, Proposition~\ref{prop:NearTwoPeripherals}
implies that
\[
   \bigcap_{A \in \Set{A}} \nbd{A}{M_2}
\]
is nonempty once $M_2$ is larger than the constant $M_1=M_1(M,\sigma)$.
Condition~(\ref{item:Unique:Center}) follows immediately.
\end{proof}

\section{Bounded packing in relatively hyperbolic groups}
\label{sec:PackingRelHyp}

Let $(G,\mathbb{P})$ be relatively hyperbolic.
An element $f \in G$ is \emph{elliptic} if it has finite order.
If $f$ has infinite order, then $f$ is \emph{parabolic} if it is contained in
a conjugate of a peripheral subgroup $P \in \mathbb{P}$
and \emph{loxodromic} otherwise.
A subgroup $H \le G$ is \emph{elliptic} if it is finite,
and \emph{parabolic} if it is infinite and contained in a conjugate of
a peripheral subgroup $P\in \mathbb{P}$.  Otherwise $H$ is \emph{loxodromic}.
Every loxodromic subgroup contains a loxodromic element.

\begin{defn}
Let $\height_{(G,\mathbb{P})}(H)=n$
if $n$ is the maximal number of distinct cosets
$g_i H$ such that the intersection $\bigcap g_i H g_i^{-1}$
is loxodromic.  Let $\width_{(G,\mathbb{P})}(H)=n$ if $n$ is the maximal number
of distinct cosets $g_i H$ such that for each $i,j$ the intersection
$g_i H g_i^{-1} \cap g_j H g_j^{-1}$ is loxodromic.
Note that an elliptic or parabolic subgroup has both height and width zero.
\end{defn}

\begin{lem}
\label{lem:ParallelLines}
Let $(G,\mathbb{P})$ be relatively hyperbolic,
and let $\Set{S}$ be a finite generating set for $G$.
There is a universal constant $L\ge 0$
so that the following holds.
Let $c$ and $c'$ be relative geodesic lines
in $\Cayley(G,\Set{S} \cup \Set{P})$
with the same endpoints at infinity.
Then the set of vertices of $c$ and the set of vertices of $c'$
are within a Hausdorff distance $L$ with respect to the metric $d_{\Set{S}}$.
\end{lem}

\begin{proof}
Since $c$ and $c'$ have the same endpoints at infinity, they are at a
finite Hausdorff distance, say $D$, in the relative metric
$d_{\Set{S} \cup \Set{P}}$.
Let $\nu$ be the constant given by Theorem~\ref{thm:ThinTriangle}.
Choose a vertex $v$ of $c$.  We will show that the $\Set{S}$--distance from
$v$ to some vertex $v'$ of $c'$ is bounded above by a constant that does not depend on $D$.

Let $x$ and $y$ be vertices of $c$ on either side of $v$,
at an $(\Set{S} \cup \Set{P})$--distance at least $D + 2\nu$
from $v$.
Choose vertices $x'$ and $y'$ of $c'$ such that
$d_{\Set{S} \cup \Set{P}}(x,x') < D$ and
$d_{\Set{S} \cup \Set{P}}(y,y') < D$.
Let $\bar{c}$ be the subpath of $c'$ from $x'$ to $y'$,
and choose relative geodesics
$c_0$ from $x$ to $x'$ and $c_1$ from $y'$ to $y$.
By Corollary~\ref{cor:ThinNGon}, the vertex $v$ is within an
$\Set{S}$--distance $2\nu$ of a vertex $v'$ of $c_0 \cup \bar{c} \cup c_1$.
However, by the triangle inequality,
the $(\Set{S} \cup \Set{P})$--distance from $v$ to any vertex of
$c_0 \cup c_1$ is greater than $2\nu$.
Therefore $v'$ is a vertex of $\bar{c}$.
\end{proof}

\begin{lem}
\label{lem:LineCloseToSubgroup}
Let $(G,\mathbb{P})$ be relatively hyperbolic with a finite generating set
$\Set{S}$.
Suppose $H\le G$ is $\sigma$--relatively quasiconvex.
If $f$ is a loxodromic element of $H$, then there exists a
bi-infinite relative geodesic $c$ in $\Cayley(G,\Set{S} \cup \Set{P})$
such that the following hold:
  \begin{enumerate}
  \item The Hausdorff $\Set{S}$--distance between the set of vertices of $c$ 
  and the subgroup $\gen{f}$ is finite.
  \item Each vertex $v$ of $c$ lies in the $\sigma$--neighborhood of $H$.
  \end{enumerate}
\end{lem}

\begin{proof}
For each $n=1,2,3,\dots$ choose a relative geodesic $c_n$ connecting
$f^{-n}$ with $f^n$.
Since $H$ is $\sigma$--relatively quasiconvex, each vertex of $c_n$
lies in the $\sigma$--neighborhood of $H$.
To complete the proof, it suffices to show that
the sequence $\{c_n\}$ has a subsequence that converges
pointwise to a relative geodesic line $c$.
In order to ensure the existence of such a subsequence
we will show that the vertices of $c_n$ track $\Set{S}$--close to
the set $\set{f^i}{-n\le i \le n}$.

The subgroup $\gen{f}$ is
$L$--relatively quasiconvex in $G$ for some $L=L(f)$,
and $\gen{f} \cap gPg^{-1}$ is finite for all conjugates of peripheral
subgroups $gP$ by Osin \cite[Theorem~4.19]{Osin06}.
Furthermore, there exists $\epsilon=\epsilon(f)$ such that
the inclusion $\gen{f} \to \Cayley(G,\Set{S} \cup \Set{P})$ 
is an $\epsilon$--quasi-isometric
embedding with respect to the word metric on $\gen{f}$
for the standard generating set $\{f\}$
(see Osin \cite[Corollary~4.20]{Osin06}).
Let $D$ be an upper bound on the
$\Set{S}$--diameters of the finitely many finite
subgroups $\gen{f} \cap gPg^{-1}$ with $\abs{g}_{\Set{P}}<L$.

By quasiconvexity, each vertex of $c_n$
lies $\Set{S}$--near $\gen{f}$.
We need to show that if $-n\le i \le n$ the vertex
$f^i$ lies $\Set{S}$--near some vertex of $c_n$.
Indeed, quasiconvexity guarantees that some of the $f^i$
lie within an $\Set{S}$--distance $\sigma'$ of vertices of $c_n$,
so we only need to bound the size of the ``gaps'' between these elements
$f^i$.

Let us begin with a pair of vertices $v$ and $w$
in the $L$--neighborhood of $\gen{f}$
connected by an edge $e$ in $\Cayley(G,\Set{S} \cup \Set{P})$.
Translate by an element of $\gen{f}$ so that
$v$ and $w$ are within an $\Set{S}$--distance $L$ of $1$ and $f^k$
for some $k\ge 0$.
If $e$ is labelled by a generator $s \in \Set{S}$ then
$\abs{f^k}_{\Set{S}}$ is at most $2L + 1$.
On the other hand, suppose $e$ is labelled by a parabolic element
$p \in P$ for some $P\in \mathbb{P}$.
Then both endpoints of $e$ lie
in one of the finitely many peripheral cosets $gP$
with $\abs{g}_{\Set{P}} < L$.
By Lemma~\ref{lem:CloseCosets} we have
\[
   \{1,f^k\} \subseteq \bignbd{\gen{f}}{L} \cap \nbd{gP}{L}
   \subseteq \bignbd{\gen{f} \cap gPg^{-1}}{L'}
\]
for some constant $L'$.  Maximizing over the finitely many possibilities
for $gP$, we can assume that $L'$ does not depend on the choice of edge $e$.
Observe that $\abs{f^k}_{\Set{S}} < 2L' + D$.
Thus in all possible cases, we have
\[
   \abs{f^k}_{\Set{S}} < M := \max\{2L+1,2L'+D\}.
\]
It follows that:
\[
   0 \le k \le \epsilon \, \abs{f^k}_{\Set{S} \cup \Set{P}} + \epsilon
     \le \epsilon \, \abs{f^k}_{\Set{S}} + \epsilon
     \le \epsilon M + \epsilon
\]

On the other hand, whenever $0\le i\le k$ we have:
\[
   \abs{f^i}_{\Set{S}} \le i\, \abs{f}_{\Set{S}} \le k \,\abs{f}_{\Set{S}}
      \le (\epsilon M + \epsilon) \,\abs{f}_{\Set{S}}
\]
Therefore the $\Set{S}$--distance from $f^i$ to $\{v,w\}$
is bounded above by the constant:
\[
   R:= L + (\epsilon M + \epsilon) \,\abs{f}_{\Set{S}}
\]
By quasiconvexity, each vertex of $c_n$ lies within an $\Set{S}$--distance
$L$ of some element of $\gen{f}$.  Since $c_n$ is an edge path,
we can apply the above argument to the set of vertices of the path $c_n$
to see that whenever $-n\le i \le n$ the vertex $f^i$ lies within an
$\Set{S}$--distance $R$ of some vertex of $c_n$.
Since $L \le R$,
the Hausdorff distance between the vertices of $c_n$
and the set $\set{f^i}{-n\le i \le n}$ is at most $R$.
An easy diagonal argument shows that the sequence of relative
geodesics $\{c_n\}$ has a subsequence converging pointwise to a relative
geodesic line $c$.

Clearly the Hausdorff $\Set{S}$--distance between the vertices
of $c$ and the subgroup $\gen{f}$ is at most $R$.
Since the vertices of $c_n$ lie in the $\sigma$--neighborhood of $H$,
the same holds for the vertices of $c$, completing the proof of the lemma.
\end{proof}

\begin{lem}[Loxodromic virtual intersection implies close]
\label{lem:LoxodromicVirtualIntersection}
Let $(G,\mathbb{P})$ be relatively hyperbolic,
let $\Set{S}$ be a finite generating set for $G$,
and choose $\sigma \ge 0$.
There is a constant $C$ so that the following holds.
Let $H$ be a $\sigma$--relatively quasiconvex subgroup of $(G,\mathbb{P})$
containing a loxodromic element $f$.
Suppose
$\Set{A} =\set{g_i H}{i \in I}$ is a collection of distinct cosets of $H$ 
such that each conjugate $g_i H g_i^{-1}$ contains some power
$f^{k_i}$ of $f$.
Then there exists a vertex $z \in G$ such that the ball of radius $C$
in $\Cayley(G,\Set{S})$ intersects every coset $g_i H \in \Set{A}$.
\end{lem}

\begin{proof}
By Osin \cite[Corollary~4.20]{Osin06}, the map $n \mapsto f^n$ is a
quasigeodesic line in $\Cayley(G,\mathcal{S} \cup \mathcal{P})$.
Therefore the set $\set{f^n}{n\in\Z}$ has exactly two limit points
$f^{\pm\infty}$ in $\boundary\Cayley(G,\mathcal{S}\cup\mathcal{P})$.
If $g_iH \in \mathcal{A}$ then $f^{k_i} \in g_iHg_i^{-1}$ for some
$k_i \in \Z$.
Thus $f^{nk_i}g_i \in g_iH$ for all $n\in\Z$.
The sets $\set{f^{nk_i}g_i}{n\in\Z}$ and $\set{f^n}{n\in\Z}$ are at a finite
Hausdorff $\Set{S}$--distance and thus have the same limit points
$f^{\pm\infty}$.

Translating $g_iH$ by $g_i^{-1}$, we see that the subgroup
$\gen{g_i^{-1} f^{k_i} g_i}$ of $H$ has limit points $g_i^{-1}f^{\pm\infty}$.
By Lemma~\ref{lem:LineCloseToSubgroup} the points $g_i^{-1}f^{\pm\infty}$
can be connected by a relative geodesic line $c_i$ such that each vertex
of $c_i$ lies within an $\Set{S}$--distance $\sigma$ of $H$.
Translating back by $g_i$ produces a relative geodesic line
$g_ic_i$ with endpoints $f^{\pm\infty}$ such that each vertex
of $g_ic_i$ lies in the $\sigma$--neighborhood of $g_iH$.

For each $i,j \in I$ the relative geodesics $g_ic_i$ and $g_jc_j$ have the
same endpoints at infinity.
By Lemma~\ref{lem:ParallelLines} the Hausdorff $\Set{S}$--distance
between the vertex sets of $g_ic_i$ and $g_jc_j$ is at most $L$.
Therefore for each vertex $z$ of $g_ic_i$, the coset $g_jH$ intersects
the ball of radius $L+\sigma$ in $\Cayley(G,\Set{S})$
centered at $z$.
\end{proof}

The following special case of
Lemma~\ref{lem:LoxodromicVirtualIntersection}
is the main connection between height/width and packing.
Roughly speaking, a collection of conjugates with loxodromic intersection
comes mutually close together in $\Cayley(G,\Set{S})$.

\begin{cor}[Loxodromic intersection implies close]
\label{cor:LoxodromicImpliesClose}
Let $(G,\mathbb{P})$ be relatively hyperbolic,
let $\Set{S}$ be a finite generating set for $G$,
and choose $\sigma \ge 0$.
There is a constant $C$ so that the following holds.
Let $H$ be a $\sigma$--relatively quasiconvex subgroup of $(G,\mathbb{P})$.
Suppose $\set{g_i H}{i \in I}$ is a collection of distinct cosets of $H$ such that
$\bigcap g_i H g_i^{-1}$ is loxodromic.
Then there exists a vertex $z \in G$ such that the ball of radius $C$
in $\Cayley(G,\Set{S})$ intersects every coset $g_i H$. \qed
\end{cor}

The following corollary is an immediate consequence of the preceding result
together with the local finiteness of $\Cayley(G,\Set{S})$.

\begin{cor}\label{cor:FiniteRelativeHeight}
Let $(G,\mathbb{P})$ be relatively hyperbolic and let $H\le G$ be
relatively quasiconvex.
Then $\height_{(G,\mathbb{P})}(H)$ is finite.\qed
\end{cor}

Corollary~\ref{cor:LoxodromicImpliesClose}
also has the following easy consequence.

\begin{cor}
\label{cor:finitely many nearby cosets}
Let $(G,\mathbb{P})$ be relatively hyperbolic and
let $H\le G$ be a relatively quasiconvex subgroup.
There are only finitely many double cosets $Hg_1H, \dots, Hg_nH$ such that
$H \cap g_i H g_i^{-1}$ is loxodromic.
\end{cor}

\begin{proof}
Suppose $\set{g_i H}{i\in I}$ is a collection of cosets such that
$H \cap g_i H g_i^{-1}$ is loxodromic for each $i$.
By Corollary~\ref{cor:LoxodromicImpliesClose},
the distances $d_{\Set{S}}(H,g_iH)$
are bounded above by a constant $C$.
Thus we can translate $g_iH$ by an element of $H$ to obtain a coset $hg_i H$
intersecting the ball of radius $C$ in $\Cayley(G,\Set{S})$ centered at the identity.
Since this ball is finite, it follows that the cosets $g_i H$ lie in only finitely many
double cosets $H g_i H$.
\end{proof}

Theorem~\ref{thm:Commensurator} also follows from
Lemma~\ref{lem:LoxodromicVirtualIntersection}.

\begin{proof}[Proof of Theorem~\ref{thm:Commensurator}]
Let $\mathcal{A}= \set{g_iH}{i\in I}$
be the set of all cosets $g_iH$ such that
$g_iHg_i^{-1} \cap H$ has finite index in both $g_iHg_i^{-1}$ and $H$.
Then the union of all cosets in $\mathcal{A}$
is the commensurator of $H$ in $G$,
and the cardinality of $\mathcal{A}$ is equal to the index of $H$
in the commensurator.
Since $H$ is loxodromic, it contains a loxodromic element $f$.
If $g_iH \in \Set{A}$ then $f^{k_i} \in g_i H g_i^{-1}$ for some
$k_i \in \Z$.
Thus by Lemma~\ref{lem:LoxodromicVirtualIntersection}
there is a vertex $z$ such that the ball of radius $C$
in $\Cayley(G,\Set{S})$ centered at $z$ intersects each coset
$g_iH \in \Set{A}$.
The finite cardinality of this ball is an upper bound for the
size of $\mathcal{A}$.
\end{proof}

\begin{prop}[Packing peripheral subgroups]
\label{prop:PeripheralPacking}
If $(G,\mathbb{P})$ is relatively hyperbolic then
each peripheral subgroup $P \in \mathbb{P}$ has bounded packing in~$G$.
\end{prop}

\begin{proof}
Choose a finite generating set $\Set{S}$ for $G$.
Let $\Set{A}$ be a collection of left cosets of $P$
whose pairwise distances are less than $D$.
Consider three distinct elements $g_0 P$, $g_1 P$, $g_2 P$ of $\Set{A}$.
Choose points $x_0 \in \nbd{g_1 P}{D} \cap \nbd{g_2 P}{D}$
and $x_1 \in \nbd{g_0 P}{D} \cap \nbd{g_2 P}{D}$.
Then by Lemma~\ref{lem:OneSideShort}, there is a constant $\zeta = \zeta(D)$
such that $d(x_0,x_1) < \zeta$.
Therefore $g_0 P$ intersects the ball of radius $\zeta + D$ centered at $x_0$.
Keeping $g_1P$ and $g_2P$ fixed and
varying $g_0 P$ over all elements of $\Set{A} - \{g_1P,g_2P\}$,
we see that every element of $\Set{A}$ intersects this ball.
Since the metric $d_{\Set{S}}$ is proper and left invariant,
the size of $\Set{A}$ is bounded in terms of $D$.
\end{proof}

\begin{thm}[Packing relative to peripheral subgroups]
\label{thm:RelPacking}
Let $(G,\mathbb{P})$ be relatively hyperbolic,
choose a finite generating set $\Set{S}$ for $G$,
and let $\Set{H}$ be a finite collection of
$\sigma$--relatively quasiconvex subgroups of $G$.
For each constant $D$, there are constants $R$ and $M_3$ such that 
the following holds.
Let $\Set{A}$ be any set of left cosets $gH$
with $g\in G$ and $H \in \Set{H}$ such that
for all $gH,g'H' \in \Set{A}$ we have $d_{\Set{S}}(gH,g'H') <D$.
Suppose the following intersection is empty:
\begin{equation}
\label{eqn:Empty}
   \bigcap_{gH \in \Set{A}} \nbd{gH}{M_3} = \emptyset
\end{equation}
Then we have the following:
\begin{enumerate}
\item \label{item:RP:Unique}
There is a unique peripheral coset $gP$ such that for all $aH \in \Set{A}$
the intersection $\nbd{gP}{R} \cap \nbd{aH}{R}$ is nonempty.
\item \label{item:RP:Triple}
$\nbd{gP}{R} \cap \nbd{aH}{R} \cap \nbd{a'H'}{R}$ is nonempty
for all $aH,a'H' \in \Set{A}$.
\item \label{item:RP:UnboundedCosets}
$\nbd{gP}{R} \cap \nbd{aH}{R}$ is unbounded
for all $aH \in \Set{A}$.
\item \label{item:RP:UnboundedConjugates}
$gPg^{-1} \cap aHa^{-1}$
is infinite for all $aH \in \Set{A}$.
\end{enumerate}
\end{thm}

\begin{proof}
We will first prove assertions (\ref{item:RP:Unique})
and~(\ref{item:RP:Triple}).
Let
\[
   n = 1 + \sum_{H \in \Set{H}} \height_{(G,\mathbb{P})}(H),
\]
which is known to be finite by Corollary~\ref{cor:FiniteRelativeHeight}.
Let $M_0 = M_0(\sigma,D,n+1)$ be the constant given by
Proposition~\ref{prop:CommonIntersection},
and let $M_2=M_2(\sigma,D,n+1,M_0)$ be the constant given by
Corollary~\ref{cor:UniquePeripheral}.
%
If there exist $n+1$ elements of $\Set{A}$
whose $M_2$--neighborhoods have empty intersection,
then Corollary~\ref{cor:UniquePeripheral}
gives (\ref{item:RP:Unique}) and~(\ref{item:RP:Triple}) immediately,
using $R=M_0$.
Thus it suffices to assume that for all
$g_1 H_1, \dots, g_{n+1} H_{n+1} \in \Set{A}$
the following intersection is nonempty:
\[
   \nbd{g_1 H_1}{M_2} \cap \dots \cap \nbd{g_{n+1} H_{n+1}}{M_2}
\]

Fix $n$ distinct cosets $g_1 H_1,\dots,g_{n} H_{n} \in \Set{A}$.
(If we choose $M_3 \ge M_2$, then (\ref{eqn:Empty}) implies that
$\Set{A}$ has at least $n+2$ elements.)
Consider the following intersection of conjugates:
\[
   K := g_1 H_1 g_1^{-1} \cap \dots \cap g_{n} H_{n} g_{n}^{-1}
\]
By our choice of $n$, some $H \in \Set{H}$ appears at least
$1+\height(H)$ times in the list $H_1,\dots,H_n$.
Thus $K$ is contained in the intersection of at least
$1+\height(H)$ distinct conjugates of $H$, so that $K$ is either
elliptic or parabolic.

\emph{Case~1:} Suppose $K$ is elliptic.
Once $M_3$ is chosen sufficiently large, we will reach a contradiction
with (\ref{eqn:Empty}).
Translating $\Set{A}$ by an element of $G$, we can assume,
without loss of generality, that the cosets $g_1 H_1, \dots, g_{n} H_{n}$
all intersect the ball of radius $M_2$ centered at the identity.
Since there are only finitely many possibilities for these cosets,
our choice of $M_3$ is permitted to depend on $g_1 H_1,\dots,g_n H_n$.
(Maximizing over all such cosets would give a uniform
constant $M_3$.)

%
%
By Corollary~\ref{cor:ManyCloseCosets} we have
\begin{equation}
\label{eqn:NearK}
   \nbd{g_1 H_1}{M_2} \cap \dots \cap \nbd{g_{n} H_{n}}{M_2}
   \subseteq \nbd{K}{M'}
\end{equation}
for some $M'$ depending on $M_2$ and the cosets
$g_1H_1,\dots,g_n H_n$.
If $aH$ is any left coset in
$\Set{A}$, the following intersection is nonempty by hypothesis:
\[
   \nbd{aH}{M_2} \cap \bigl( \nbd{g_1 H_1}{M_2} \cap \dots
                           \cap \nbd{g_{n} H_{n}}{M_2} \bigr).
\]
Thus by (\ref{eqn:NearK}) we have $d_{\Set{S}}(aH,K) < M_2 + M'$.
Since $K$ is an elliptic subgroup,
it has a finite diameter $D=D(K)$ and contains the identity.
Hence every $aH \in \Set{A}$ intersects the finite ball of radius
$D+M_2+M'$ centered at the identity,
which contradicts (\ref{eqn:Empty}) if we choose $M_3 \ge D+M_2+M'$.



\emph{Case~2:}
Suppose $K$ is parabolic.
By Proposition~\ref{prop:FinConjClasses}, for each $H \in \Set{S}$
and each peripheral subgroup $P$,
the peripheral cosets $gP$ such that $H \cap gPg^{-1}$ is infinite
lie in finitely many $H$--orbits.
Let $T$ be the maximum distance between $H$ and any coset $gP$
such that $H \cap gPg^{-1}$ is infinite.
Equivalently, if $gHg^{-1} \cap P$ is infinite, then
the distance between $gH$ and $P$ is at most $T$.

Since $K$ is parabolic, it is contained in a conjugate of some
peripheral subgroup $P$.  Translating $\Set{A}$ by an element of $G$,
we can assume that $K$ is a subgroup of $P$ itself.
For each $i=1,\dots,n$ the group $g_i H_i g_i^{-1} \cap P$
is infinite, since it contains the infinite group $K$.
Therefore the distance between $g_i H_i$ and $P$ is at most $T$.
By assumption, the following intersection is nonempty:
\[
   \nbd{g_1 H_1}{M_2} \cap \dots \cap \nbd{g_n H_n}{M_2}
\]
Thus by Corollary~\ref{cor:AllNearPeripheral}
there is a constant $R_1=R_1(\sigma,T+M_2,n)$ and a point
\[
   x \in \nbd{P}{R_1} \cap \nbd{g_1 H_1}{R_1} \cap \cdots
      \cap \nbd{g_n H_n}{R_1}.
\]
Translating $\Set{A}$ by an element of $P$
we can assume without loss of generality that
$x$ lies in the ball of radius $R_1$ centered at the identity.
Consequently, the cosets $g_1 H_1,\dots,g_n H_n$ each intersect the
ball of radius $2R_1$ centered at the identity.
As in the elliptic case, there are a uniformly bounded number of choices
for $g_1 H_1,\dots,g_n H_n$.

As before, (\ref{eqn:NearK}) holds for some $M'$
depending on $g_1 H_1,\dots,g_{n}H_{n}$ and $M_2$.
Since $K\le P$, it follows from the same reasoning as above that
$d_{\Set{S}}(aH,P) \le d_{\Set{S}}(aH,K) < M_2+M'$ for each $aH \in \Set{A}$.
Furthermore, for each $aH,a'H' \in \Set{A}$, the pairwise distances
between $aH$, $a'H'$ and $P$ are at most $D+M_2+M'$.
Thus Proposition~\ref{prop:Center}(\ref{item:Center2a})
gives a constant $R=R(\sigma,D+M_2+M')$ such that
for all $aH,a'H' \in \Set{A}$
the following intersection is nonempty:
\[
   \nbd{P}{R} \cap \nbd{aH}{R} \cap \nbd{a'H'}{R}
\]
Recall that $R$ depends on the cosets $g_1H_1,\dots,g_n H_n$,
which each
intersect the ball of radius $2R_1$ centered at the identity.
Maximizing over the finitely many possible choices for these cosets
gives a uniform constant---which we again call $R$---that does not depend
on our particular choice of $\Set{A}$.

In order to prove (\ref{item:RP:Unique}) and~(\ref{item:RP:Triple})
it suffices to establish the uniqueness of $P$ in the sense of
(\ref{item:RP:Unique}) once $M_3$ is chosen sufficiently large.
Suppose by way of contradiction
that there were another peripheral coset $g'P' \ne P$
such that for all $aH \in \Set{A}$ the intersection
$\nbd{g'P'}{R} \cap \nbd{aH}{R}$ is nonempty.
Then Proposition~\ref{prop:NearTwoPeripherals}
implies that the intersection
$\bigcap_{aH\in\Set{A}} \nbd{aH}{M_1}$
is nonempty for $M_1=M_1(R,\sigma)$.
If we choose $M_3$ larger than $M_1$, 
this clearly contradicts (\ref{eqn:Empty}).


We will now show that (\ref{item:RP:Triple})
implies (\ref{item:RP:UnboundedCosets}) when $M_3$ is chosen
sufficiently large.
Indeed, suppose by way of contradiction
that there exists $aH \in \Set{A}$
such that $\nbd{gP}{R} \cap \nbd{aH}{R}$ is bounded.
After translating $\Set{A}$ by an element of $G$, we can assume
that $gP$ and $aH$ both intersect the ball of radius $R$ in $G$.
Thus there are a bounded number of choices for the cosets $gP$
and $aH$, up to translation in $G$.
Maximizing over all such choices gives a uniform upper bound
$S$ on the diameter of the intersection $I:=\nbd{gP}{R} \cap \nbd{aH}{R}$.
If $a'H'$ is any other element of $\Set{A}$, then
(\ref{item:RP:Triple}) implies that $\nbd{a'H'}{R}$ intersects $I$.
Therefore every $a'H' \in \Set{A}$ intersects a fixed set of diameter $S+R$,
which contradicts (\ref{eqn:Empty}) if we choose $M_3$
larger than $S+R$.

Finally observe that (\ref{item:RP:UnboundedConjugates})
follows immediately from (\ref{item:RP:UnboundedCosets})
using Lemma~\ref{lem:CloseCosets}.
\end{proof}

\begin{thm}[Packing relatively quasiconvex subgroups]
\label{thm:PackingRelQC}
Let $(G,\mathbb{P})$ be relatively hyperbolic,
choose a finite generating set $\Set{S}$ for $G$,
and let $H$ be a $\sigma$--relatively quasiconvex subgroup of $G$.
Suppose for each peripheral subgroup $P \in \mathbb{P}$ and each $g \in G$
the intersection $P \cap gHg^{-1}$ has bounded packing in $P$.
Then $H$ has bounded packing in $G$.
\end{thm}

\begin{proof}
Let $\Set{A}$ be any set of left cosets $gH$ whose pairwise distances
are at most $D$.
Our goal is to show that $\Set{A}$ is finite and to bound the cardinality
of $\Set{A}$ as a function of $D$.

Let $R$ and $M_3$ be the constants given by Theorem~\ref{thm:RelPacking}.
Suppose the following intersection is nonempty:
\[
   \bigcap_{gH\in\Set{A}} \nbd{gH}{M_3}
\]
Then there is a fixed ball of radius $M_3$ in $G$
that intersects every $gH \in \Set{A}$.
The cardinality of such a ball is an upper bound on the cardinality
of $\Set{A}$.

Thus it suffices to assume that the following intersection is empty:
\[
   \bigcap_{gH\in\Set{A}} \nbd{gH}{M_3}
\]
In this case, Theorem~\ref{thm:RelPacking} applies.
After translating
$\Set{A}$ by an element of $G$, there exists a
peripheral subgroup $P$ such that for all $gH,g'H \in \Set{A}$
the following intersection is nonempty:
\[
   \nbd{P}{R} \cap \nbd{gH}{R} \cap \nbd{g'H}{R}
\]

For each $gH \in \Set{A}$, there is a point $p \in P$ such that
$d_{\Set{S}}(p,gH) < 2R$.
Translating by $p^{-1}$, we see that the coset $p^{-1}gH$
intersects the ball of radius $2R$ centered at the identity.
Therefore every coset in $\Set{A}$ can be expressed as $pgH$
for some $p \in P$ and $\abs{g}_{\Set{S}} < 2R$.
The cosets $gH \in \Set{A}$
lie in a bounded number of orbits under the action of $P$;
therefore it is enough to bound the number of elements of $\Set{A}$
lying in a single orbit $P(gH)$; in other words, a double coset.

Setting $Q := P \cap gHg^{-1}$ the map
$pgH \mapsto pQ$ that takes a left coset of $H$ in $P(gH)$ to the corresponding coset of $Q$ in $P$ is well-defined and injective;
for $pgH = p'gH$ if and only if $p^{-1}p' \in P \cap gHg^{-1}=Q$
if and only if $pQ=p'Q$.
By Lemma~\ref{lem:CloseCosets}, we have
\[
   \nbd{P}{R} \cap \nbd{gH}{R} \subseteq \nbd{Q}{M'}
\]
for some $M'$ depending on $P$, $gH$ and $R$.
Translating by $p$ gives
\[
   \nbd{P}{R} \cap \nbd{pgH}{R} \subseteq \nbd{pQ}{M'}
\]
for each $p \in P$, as illustrated in Figure~\ref{fig:footprints}.

\begin{figure}
\labellist
\small\hair 2pt
\pinlabel $P$ at 125 8
\pinlabel $pQ$ at 71 28
\pinlabel $p'Q$ at 164 30
\pinlabel $pgH$ at 74 95
\pinlabel $p'gH$ at 164 101
\endlabellist
\begin{center}
\includegraphics{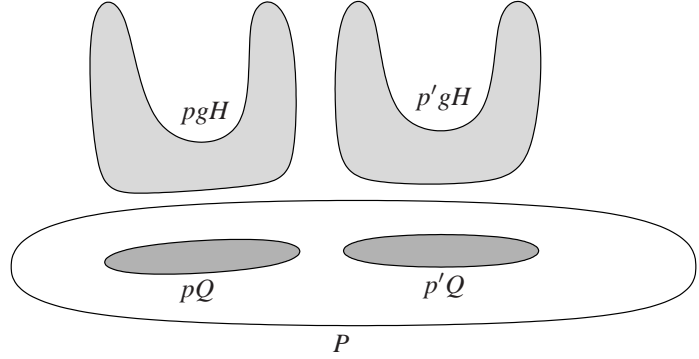}
\end{center}
\caption{The coset $pQ$ is coarsely the intersection of
neighborhoods of $P$ and $pgH$.  The existence of a point close to
$P$, $pgH$ and $p'gH$ implies that $pQ$ and $p'Q$ are close.}
\label{fig:footprints}
\end{figure}

By hypothesis, for each $pgH, p'gH \in \Set{A}$ we have
\begin{multline*}
   \nbd{P}{R} \cap \nbd{pgH}{R} \cap \nbd{p'gH}{R} \\
     = \bigl( \nbd{P}{R} \cap \nbd{pgH}{R} \bigr)
       \cap \bigl( \nbd{P}{R} \cap \nbd{p'gH}{R} \bigr) \ne \emptyset.
\end{multline*}
Therefore
\[
   \nbd{pQ}{M'} \cap \nbd{p'Q}{M'} \ne \emptyset.
\]
But $Q=P\cap gHg^{-1}$ has bounded packing in $P$
so we have an upper bound on the number of cosets $pQ$ in $P$
that correspond to cosets $pgH \in \Set{A}$.
\end{proof}

\begin{cor}\label{cor:RelHypFiniteWidth}
Let $(G,\mathbb{P})$ be relatively hyperbolic and $H\le G$ a relatively 
quasiconvex subgroup.
Suppose $P \cap gHg^{-1}$ has bounded packing in $P$ for each peripheral subgroup $P \in \mathbb{P}$ and
each $g \in G$.
Then $\width_{(G,\mathbb{P})}(H)$ is finite.
\end{cor}

\begin{proof}
If $g_0 H$ and $g_1 H$ are distinct cosets such that
$g_0 H g_0^{-1} \cap g_1 H g_1^{-1}$
is loxodromic, then Corollary~\ref{cor:LoxodromicImpliesClose}
implies that $d_{\mathcal{S}}(g_0 H, g_1 H) < 2C$.
The result follows immediately since $H$ has bounded packing in $G$
by Theorem~\ref{thm:PackingRelQC}.
\end{proof}

\section{Quasiconvex and separable implies virtually malnormal}
\label{sec:SeparableMalnormal}

In this section we present a malnormality consequence of
separability and finite height. The results are new even in the
word-hyperbolic case.

\begin{defn}
A subgroup $H$ of $G$ is \emph{malnormal} if $H \cap gHg^{-1}$ is trivial
for each $g\in G-H$.
In applications of malnormality it is usually sufficient to know
that the intersection of conjugates is finite instead of trivial,
and so we have elsewhere defined $H$ to be \emph{almost malnormal}
if $H\cap gHg^{-1}$ is finite for each $g\in G-H$.

In light of our earlier definitions of width and height of subgroups
of a relatively hyperbolic group $G$,
it is natural to define a subgroup $H$ of $G$
to be \emph{relatively malnormal} if
$H\cap gHg^{-1}$ is either elliptic or parabolic for each $g\in G-H$.

Note that any malnormal subgroup and any parabolic or elliptic subgroup is immediately
relatively malnormal.
\end{defn}

\begin{defn}
A subgroup $H$ of a group $G$ is \emph{separable}
if $H$ is an intersection of finite index subgroups of $G$.
\end{defn}

\begin{thm}
\label{thm:RelMalnormal}
Let $H$ be a separable, relatively quasiconvex subgroup of the relatively 
hyperbolic group $G$.
Then there is a finite index subgroup $K$ of $G$ containing $H$
such that $H$ is relatively malnormal in $K$.
\end{thm}

\begin{proof}
By Corollary~\ref{cor:finitely many nearby cosets},
there are only finitely many double cosets $Hg_1H, \dots, Hg_nH$ such that
$H\cap g_i H g_i^{-1}$ is loxodromic.
Since $H$ is separable, there exists a finite index subgroup $K$ of $G$
such that $H \le K$ but $g_i\not\in K$ for $1\leq i\leq n$.
If $k \in K-H$ and $H \cap kHk^{-1}$ is loxodromic, then
$kH = hg_iH$ for some $g_i$ and some $h \in H$.
Therefore $hg_i H \subset K$ and $H\le K$ so $g_i \in K$, contradicting our choice of $K$.
Consequently $H$ is relatively malnormal in $K$.
\end{proof}

\begin{prob}\label{prob:malnormal}
Find a relatively hyperbolic group $G$
with a relatively malnormal subgroup $M$
that is not quasiconvex.
\end{prob}

There is a related long-standing problem of constructing
 a malnormal subgroup of a word-hyperbolic
group that is not quasiconvex. One expects that Problem~\ref{prob:malnormal}
won't be much easier, even if we relax the condition of malnormality
to merely require finite width.

\bibliographystyle{gtart}
\bibliography{Packing}

\def\polhk#1{\setbox0=\hbox{#1}{\ooalign{\hidewidth
  \lower1ex\hbox{$\,\lhook$}\hidewidth\crcr\unhbox0}}}
  \def\RomanianComma#1{\setbox0=\hbox{#1}{\ooalign{\hidewidth
  \lower1.2ex\hbox{$\mspace{1mu}^{,}$}\hidewidth\crcr\unhbox0}}}
\begin{thebibliography}{}
\providecommand\bibmarginpar{\leavevmode\marginpar}
\def\urlstyle#1{{\tt #1}}

\bibitem{ABC91}
\textbf{J Alonso}, \textbf{T Brady}, \textbf{D Cooper}, \textbf{V Ferlini},
  \textbf{M Lustig}, \textbf{M Mihalik}, \textbf{H Short}, \emph{Notes on word
  hyperbolic groups}, (H Short, editor), from: ``Group theory from a geometrical viewpoint
  \textup{(}Trieste, 1990\textup{)}'', ({\'E} Ghys, A Haefliger, A Verjovsky,
  editors), World Sci.\ Publishing, River Edge, NJ (1991)  3--63

\bibitem{BandeltVanDeVel91}
\textbf{H-J Bandelt}, \textbf{M van~de Vel}, \emph{Superextensions and the
  depth of median graphs}, J. Combin. Theory Ser. A 57 (1991) 187--202

\bibitem{BowditchRelHyp}
\textbf{B Bowditch}, \emph{Relatively hyperbolic groups} (1999) Preprint, Univ.\
  of Southampton

\bibitem{Chepoi00}
\textbf{V Chepoi}, \emph{Graphs of some {$\CAT(0)$} complexes}, Adv. in Appl.
  Math. 24 (2000) 125--179

\bibitem{DrutuSapirTreeGraded}
\textbf{C Dru{\RomanianComma{t}}u}, \textbf{M Sapir}, \emph{Tree-graded spaces
  and asymptotic cones of groups}, with an appendix
  by D~Osin and M~Sapir,
Topology 44 (2005) 959--1058

\bibitem{Farb98}
\textbf{B Farb}, \emph{Relatively hyperbolic groups}, Geom.\ Funct.\ Anal. 8
  (1998) 810--840

\bibitem{Gerasimov97}
\textbf{V Gerasimov}, \emph{Semi-splittings of groups and actions on cubings},
  from: ``Algebra, geometry, analysis and mathematical physics
  \textup{(}Russian\textup{)} \textup{(}Novosibirsk, 1996\textup{)}'', Izdat.
  Ross. Akad. Nauk Sib. Otd. Inst. Mat., Novosibirsk (1997)  91--109 English
  translation published as \emph{Fixed-point-free actions on cubings},
  Siberian Adv. Math. 8 (1998) 36--58

\bibitem{GMRS98}
\textbf{R Gitik}, \textbf{M Mitra}, \textbf{E Rips}, \textbf{M Sageev},
  \emph{Widths of subgroups}, Trans. Amer. Math. Soc. 350 (1998) 321--329

\bibitem{Gromov87}
\textbf{M Gromov}, \emph{Hyperbolic groups}, from: ``Essays in group theory'',
  (S Gersten, editor), Springer, New York (1987)  75--263

\bibitem{HruskaRelQC}
\textbf{G Hruska}, \emph{Relative hyperbolicity and relative quasiconvexity for
  countable groups} arXiv:0801.4596v1 [math.GR]

\bibitem{HruskaWiseAxioms}
\textbf{G Hruska}, \textbf{D Wise}, \emph{Finiteness properties of cubulated
  groups} In preparation

\bibitem{MartinezPedroza}
\textbf{E Martinez-Pedroza}, \emph{Combination of quasiconvex subgroups of
  relatively hyperbolic groups} To appear in \textit{Groups Geom. Dyn.}
  arXiv:0709.2113v2 [math.GR]

\bibitem{NibloReeves03}
\textbf{G Niblo}, \textbf{L Reeves}, \emph{Coxeter groups act on {$\CAT(0)$}
  cube complexes}, J. Group Theory 6 (2003) 399--413

\bibitem{NibloRoller98}
\textbf{G Niblo}, \textbf{M Roller}, \emph{Groups acting on cubes and
  {K}azhdan's property \textup{(}{T}\textup{)}}, Proc.\ Amer.\ Math.\ Soc. 126 (1998) 693--699

\bibitem{Osin06}
\textbf{D Osin}, \emph{Relatively hyperbolic groups: {I}ntrinsic geometry,
  algebraic properties, and algorithmic problems}, Mem. Amer. Math. Soc. 179
  (2006) 1--100

\bibitem{Rips82}
\textbf{E Rips}, \emph{Subgroups of small cancellation groups}, Bull. London
  Math. Soc. 14 (1982) 45--47

\bibitem{RollerPocSets}
\textbf{M Roller}, \emph{Poc sets, median algebras and group actions. {A}n
  extended study of {D}unwoody's construction and {S}ageev's theorem}
  (1998) Preprint, Univ.\ of Southampton

\bibitem{RubinsteinSageev99}
\textbf{H Rubinstein}, \textbf{M Sageev}, \emph{Intersection patterns of
  essential surfaces in {$3$}--manifolds}, Topology 38 (1999) 1281--1291

\bibitem{RubinsteinWang98}
\textbf{J Rubinstein}, \textbf{S Wang}, \emph{{$\pi\sb 1$}--injective surfaces
  in graph manifolds}, Comment. Math. Helv. 73 (1998) 499--515

\bibitem{Sageev95}
\textbf{M Sageev}, \emph{Ends of group pairs and non-positively curved cube
  complexes}, Proc.\ London Math.\ Soc.\ \textup{(}3\textup{)} 71 (1995)
  585--617

\bibitem{Sageev97}
\textbf{M Sageev}, \emph{Codimension--$1$ subgroups and splittings of groups},
  J.\ Algebra 189 (1997) 377--389

\end{thebibliography}

\end{document}